\documentclass{article}
\usepackage{amsmath,amssymb,amsfonts}
\usepackage{graphicx}
\usepackage{overpic}
\usepackage{wrapfig}
\usepackage{color}
\usepackage{caption}
\usepackage{subcaption}

\newcommand{\bu}{\mathbf u}
\newcommand{\bn}{\mathbf n}
\newcommand{\tbn}{\mathbf {\tilde n}}
\newcommand{\bv}{\mathbf v}
\newcommand{\by}{\mathbf y}
\newcommand{\bV}{\mathbf V}
\newcommand{\wbV}{\widetilde{V}_h}
\newcommand{\bx}{\mathbf x}

\newcommand{\bI}{\mathbf I}
\newcommand{\bD}{\mathbf D}

\newcommand{\be}{\mathbf e}

\newcommand{\T}{\mathcal T}

\newcommand{\bsigma}{\boldsymbol{\sigma}}

\newcommand{\eps}{\varepsilon}

\newcommand{\bmat}[1]{\begin{bmatrix} #1 \end{bmatrix}}

\newtheorem{proposition}{Proposition}
\newtheorem{remark}{Remark}

\newcommand{\Div}{\operatorname{div}}
\def\dO{\partial\Omega }
\newcommand{\dom}{\partial\omega }

\usepackage{fullpage}
\usepackage{cite}

\usepackage{accents}
\renewcommand*{\dot}[1]{%
	\accentset{\mbox{\large\bfseries .}}{#1}}

\newcommand\rev[1]{{\color{black}#1}}

\begin{document}

\title{Local conservation laws of continuous Galerkin method for the incompressible Navier--Stokes equations in EMAC form}
\author{Maxim A. Olshanskii\thanks{\small Department of Mathematics, University of Houston, Houston, TX, 77204, maolshanskiy@uh.edu.}
\and
Leo G. Rebholz
\thanks{\small School of Mathematical and Statistical Sciences, Clemson University, Clemson, SC, 29364, rebholz@clemson.edu.}
}
\date{}

\maketitle

\begin{abstract}
We consider {\it local} balances of momentum and angular momentum for the incompressible Navier-Stokes equations.  First, we formulate new weak forms of the physical balances (conservation laws) of these quantities, and prove they are equivalent to the usual conservation law formulations.  We then show that \emph{continuous} Galerkin discretizations of the Navier-Stokes equations using the EMAC form of the nonlinearity preserve discrete analogues of the weak form conservation laws, both in the Eulerian formulation and the Lagrangian formulation (which are not equivalent after discretizations).  Numerical tests illustrate the new theory. 
\end{abstract}

\section{Introduction}

We are interested in conservation properties of continuous Galerkin discretizations of the incompressible Navier-Stokes equations (NSE), which are given by
\begin{equation}\label{NSE}
\left\{
\begin{aligned}
\frac{\partial \bu}{\partial t}+(\bu\cdot \nabla)\bu - \Div \bsigma &=0 \\
\Div \bu&=0
\end{aligned}\right.\quad\text{in}~~\Omega.
\end{equation}
Here, $\bsigma$ is the Cauchy stress tensor, and we restrict to the case of a Newtonian fluid with  $\bsigma=2\nu\bD(\bu)-p\bI$, where $\bD(\bu)=\tfrac{1}{2}(\nabla\bu+\nabla^T\bu)$ is a rate of deformation tensor.

It is well known that the smooth solution to \eqref{NSE} obeys an array of conservation laws,  including the conservation of momentum, energy, vorticity, etc., which can be expressed in terms of proper balances for material volumes of fluid. The development of numerical methods that provide discrete counterparts for possibly many of these conservation laws is a long-standing challenge for the computational fluid dynamics community. This challenge has been addressed by numerous authors and from various perspectives; for example, see \cite{A66, AM03,  LW04,OR10b, evans2013isogeometric,SCN15,PG16,CHOR17,coppola2019discrete} and references therein.
Many of these studies have considered the \emph{global} conservation properties of numerical methods, i.e., balances of physical quantities across the entire computational domain. While properly calibrating these global integral statistics is necessary for a method to be long-time accurate, it is difficult to see how this alone can guarantee the quality of a numerical solution.

The proper \emph{local} balances of momentum, energy, vorticity, etc. represent a significantly stronger requirement for a numerical solution. Note that ``element-wise conservation" is a common argument used to motivate the application of discontinuous Galerkin or finite volume discretization techniques (see, for instance, \cite{cockburn2003discontinuous,leveque2002finite}). At the same time, there is a widespread belief that continuous (velocity $H^1$-conforming) Galerkin methods inevitably violate local conservation laws; however, see \cite{HEML00,HW05} for a different viewpoint.

Another obstacle in achieving proper discrete counterparts of both local and global conservation laws for $\bu$ is the fact that continuous Galerkin discretizations of \eqref{NSE} (e.g., conforming finite element methods) typically enforce the divergence-free constraint only weakly \cite{CHOR17}. The purpose of this paper is to demonstrate that a continuous Galerkin solution, which is only weakly divergence-free, for \eqref{NSE} does  satisfy properly formulated local conservation laws for momentum and angular momentum when one applies the so-called EMAC (Energy, Momentum, and Angular Momentum Conserving) formulation of the NSE.

The EMAC formulation of the discrete NSE was originally developed in \cite{CHOR17}. It re-writes the inertia terms as
\[
\bu\cdot\nabla \bu \rightarrow 2\bD(\bu_h)\bu_h + (\Div \bu_h)\bu_h,
\]
along with an altered pressure $p_h$ representing $p-\frac12 |\bu|^2$.  The motivation for EMAC was that Galerkin schemes using it can be shown to conserve global energy, momentum and angular momentum balances when $\Div \bu_h\ne 0$, while schemes using the common nonlinearity formulations such as convective (CONV: $\bu_h\cdot\nabla \bu_h$), skew-symmetric (SKEW: $\bu_h\cdot\nabla \bu_h + \frac12(\Div \bu_h)\rev{\bu_h}$) and rotational (ROT: $(\nabla \times \bu_h) \times \bu_h$) do not preserve some or all of  these quantities.  Perhaps not surprisingly, the use of EMAC has become popular for large scale fluid computations in a wide variety of applications and is shown to give better accuracy especially over longer time intervals e.g. \cite{PCLRH18,LHOCR19,SP18,SP18b,LPH19,KMLP23,VSA22,CHOR17,CHOR19,OR20,INRRV23} and is built into Alya which is a massively parallel multiphysics unstructured finite element code~\cite{VHK16}.  In addition to the better discrete physics of EMAC discussed above, it was proven in \cite{OR20} that schemes using EMAC are more long time stable because the Gronwall constant can be shown to be independent of (explicit) dependence on the Reynolds number ($Re$), \rev{and in \cite{garcia2021convergence} the uniform in $Re$ error estimate was derived for the EMAC error}; such results are not known for skew-symmetric, convective, or rotation forms for commonly used velocity-pressure finite elements such as Taylor--Hood elements.  

The purpose of this paper is to provide more theoretical justification that EMAC is superior compared to other discrete nonlinearity formulations, by proving that 
continuous Galerkin discretizations using EMAC admit an exact {\it local} balance of momentum and angular momentum.  There are very few results for local conservation properties of continuous   finite element methods, with \cite{HEML00,HW05} being two fundamental works in this direction.  The paper \cite{HW05} showed that for NSE, typical Galerkin schemes are not generally conservative, although this can be `fixed' by multiscale formulation and adding a residual term.
One observation made in this paper is that although local balances written in different forms --   standard Eulerian, Lagrangian, or weak Eulerian and Lagrangian forms introduced here -- 
represent the same conservation laws of fluid momentum and angular momentum,  after discretization each form can be different.  By considering the weak forms, which we refer to as \emph{diffuse-volume forms}, of conservation laws, we can demonstrate that EMAC continuous Galerkin discretizations exactly preserve properly formulated local momentum and angular momentum balances.  Furthermore, the discrete balances established here serve as direct analogies to the balances at the partial differential equation (PDE) level, obviating the need for a multiscale approach and additional residual terms to establish this connection.  We note also that from the proof construction, it is not possible for SKEW, CONV or ROT to preserve these local balances of momentum and angular momentum in the same manner that EMAC does, since they do not preserve them globally.

The rest of the paper is arranged as follows. Section~\ref{S2} recalls local conservation laws of momentum and angular momentum. 
The laws can be equivalently formulated in Eulerian and Lagrangian forms.  Section~\ref{S3} introduces a different way to formulate local conservation laws, which is given the name {diffuse-volume} form of the   conservation laws due to some similarity with diffuse-interface or phase-field methods in fluid mechanics. We show that this is another equivalent way to formulate the local balances.  Here we also distinguish between diffuse-volume Eulerian and Lagrangian forms.
 Section~\ref{S4} demonstrates how the continuous Galerkin method for the NSE in EMAC form satisfies discrete counterparts of the (local)   diffuse-volume Eulerian and Lagrangian conservation laws. Section~\ref{S5} offers a few illustrative numerical examples.

\section{Eulerian and Lagrangian forms of momentum and angular momentum conservation}\label{S2}

To formulate {local} conservation laws satisfied by a \emph{smooth} solution to \eqref{NSE}, we fix some $t$ and let $\omega\subset\Omega$ be a \emph{fixed} subdomain of $\Omega$ with sufficiently smooth  boundary $\dom$. \rev{We shall assume that $\Omega$ is bounded}.
For this volume $\omega$, the   balance of  momentum and angular momentum take the form:
	\begin{align}
		\text{Moment.}\quad & \frac{d}{dt}\int_{\omega}\bu\,dx&=&~ 2\nu\int_{\dom}\bD(\bu)\bn\,ds&-&\quad\int_{\dom}p\bn\,ds&-&\int_{\dom}\bu(\bu\cdot\bn)\,ds, 
		\label{Laws2m}\\ 
		\text{Ang. Moment.}\quad & \underbrace{\frac{d}{dt}\int_{\omega}\bu\times\bx\,dx}_{\scriptsize\begin{array}{c}\text{momentum} \\ \text{rate of change}\end{array}}& =&~\underbrace{2\nu\int_{\dom}(\bD(\bu)\bn)\times\bx\,ds}_{\scriptsize\begin{array}{c}\text{momentum change due}\\ \text{to friction on }\dom\end{array}}&-&\underbrace{\int_{\dom}p(\bn\times\bx)\,ds}_{\scriptsize\begin{array}{c}\text{moment. change due}\\ \text{to pressure on }\dom\end{array}}&-&~\underbrace{\int_{\dom}(\bu\times\bx)(\bu\cdot\bn)\,ds.}_{\scriptsize\begin{array}{c}\text{the flux of}\\ \text{momentum through }\dom\end{array}} \label{Laws2a}
	\end{align}
Hereafter $\bn$ is the outward normal vector on $\dom$.

Local balances~\eqref{Laws2m}--\eqref{Laws2a} can be interpreted as \emph{Eulerian form} of the conservation laws, in contrast to the \emph{Lagrangian form} formulated for a material volume below.
\smallskip

We now let $\Omega_t\subset\Omega$ be a \emph{material} volume of the fluid. 
For the material volume, the   conservation laws for momentum and angular momentum take the form:
	\begin{align}
		\text{Momentum}\quad & \frac{d}{dt}\int_{\Omega_t}\bu\,dx&=&~  2\nu\int_{\dO_t}\bD(\bu)\bn\,ds-\int_{\dO_t}p\bn\,ds,& 
		\label{Laws1m}\\
\hspace{5ex}		\text{Angular Momentum}\quad & \frac{d}{dt}\int_{\Omega_t}\bu\times\bx\,dx&=&~ 2\nu\int_{\dO_t}(\bD(\bu)\bn)\times\bx\,ds-\int_{\dO_t}p(\bn\times\bx)\,ds.& \hspace{5ex} \label{Laws1a}
	\end{align}

Of course, for smooth solutions to \eqref{NSE} the Eulerian and Lagrangian forms are just two different formulations of the same fundamental laws of continuum mechanics. They both follow from \eqref{NSE}, and conversely, together with mass conservation they imply \eqref{NSE}.  This  equivalence of \eqref{NSE} to the validity of local conservation laws (specifically, those concerning mass and momentum) is textbook material.  The standard tools used to verify this equivalence include the divergence theorem, the freedom to choose fluid volumes $\Omega_t$ or $\omega$,
and the Reynolds' transport theorem to handle the Lagrangian form, which states
\[
\frac{d}{dt}\int_{\Omega_t}f\,dx = \int_{\Omega_t}\Big(\frac{\partial f}{\partial t}+\Div(f\bu)\Big)\,dx,
\] 
for a smooth scalar function $f$.
\smallskip

Continuous Galerkin  methods like the $H^1$-conforming  finite element method (FEM) employ finite dimensional  subspaces of Sobolev spaces to project \eqref{NSE} and typically do not offer enough flexibility to verify a direct analogue of \eqref{Laws2m}--\eqref{Laws2a} or \eqref{Laws1m}--\eqref{Laws1a}. Below we reformulate local conservation laws in a form more convenient for continuous Galerkin methods.

\smallskip

\section{Weak form of the conservation laws}\label{S3}
The purpose of this section is to derive the conservation laws \eqref{Laws2m}--\eqref{Laws2a} and \eqref{Laws1m}--\eqref{Laws1a} in a form {more appropriate} for a variational formulation. Let $\omega\subset\Omega$ be an arbitrary subdomain of $\Omega$ with sufficiently smooth $\dom$.
Denote by $\phi$ an \emph{arbitrary} smooth  function  such that $\omega=\mbox{supp}(\phi)$, and set
\[
\tbn:=-\frac{\nabla\phi}{|\nabla\phi|}
\]
for $\bx$ such that $\nabla\phi(\bx)\neq0$, and let $\tbn(\bx)$ be an arbitrary vector of unit length if  $\nabla\phi(\bx)=0$. Note that 
$\tbn(\bx)=\bn(\bx)$ for $\bx\in\dom$. 
To obtain the weak  form of the laws, we multiply the first equation in \eqref{NSE} by $\phi\be_i$ for  momentum conservation and by 
$\phi\be_i\times\bx$ for angular  momentum conservation. Doing this for $i=1,\dots,d$, integrating over $\omega$ and by parts leads after 
straightforward computations to the following \emph{weak form of the conservation laws}:
	\begin{align}
		\text{Moment.}\quad & \frac{d}{dt}\int_{\omega}\phi\bu\,dx\hskip4ex = 2\nu\int_{\omega}\bD(\bu)\tbn|\nabla\phi|\,dx-\int_{\omega}p\tbn|\nabla\phi|\,dx-\int_{\omega}\bu(\bu\cdot\tbn)|\nabla\phi|\,dx, \label{Laws2phim}\\
		\text{Ang. Moment.}\quad & \frac{d}{dt}\int_{\omega}\phi\bu\times\bx\,dx =2\nu\int_{\omega}(\bD(\bu)\tbn)\times\bx\,|\nabla\phi|\,dx-\int_{\omega}p(\tbn\times\bx)|\nabla\phi|\,dx-\int_{\omega}(\bu\times\bx)(\bu\cdot\tbn)|\nabla\phi|\,dx.  \label{Laws2phia}
	\end{align}
We note that for all calculations below to make sense it is sufficient to assume $\phi\in W^{1,\infty}(\Omega)$. 

\medskip 
Given the freedom in choosing  $\omega$ and $\phi$ one can show that \eqref{Laws2m}--\eqref{Laws2a} and \eqref{Laws2phim}--\eqref{Laws2phia} are \emph{equivalent} if $\bu$ is sufficiently smooth and divergence free. We formulate it as a proposition.

\begin{proposition}\label{pr1}
Assume $\bu$ and $p$ are smooth and $\Div\bu=0$, then \eqref{Laws2m} {\rm (or \eqref{Laws2a})} holds for any subdomain $\omega\subset\Omega$ iff  \eqref{Laws2phim} {\rm (or \eqref{Laws2phia})} holds for any $\phi\in W^{1,\infty}(\Omega)$ with $\mbox{supp}(\phi)\subset\Omega$.
\end{proposition}
\noindent\textit{Proof.}
We know that \eqref{Laws2m} and $\Div\bu=0$ imply \eqref{NSE} by standard arguments, given that $\omega$ can be taken as an arbitrary subdomain of $\Omega$ and for any $t$. Similarly, the fact that \eqref{Laws2phim} holds for any $\phi\in \dot{C}(\Omega)$ leads to 
\[
\int_\Omega\Big(\frac{\partial \bu}{\partial t}+(\bu\cdot \nabla)\bu - 2\nu\Div\bD(\bu) +\nabla p \Big)\phi\,dx=0,\quad\forall\, \phi\in \dot{C}(\Omega),
\] 
which implies \eqref{NSE} due to the density of smooth compactly supported functions in $L^2(\Omega)$.   In turn, both \eqref{Laws2m} and \eqref{Laws2phim} are quick consequences of \eqref{NSE}. Thus \eqref{Laws2m} implies  \eqref{Laws2phim} and vice versa.

The same arguments can be applied to show the equivalence of \eqref{Laws2a} and \eqref{Laws2phia}. The only difference is that the equivalence of  \eqref{Laws2a} is established not to the momentum equation in \eqref{NSE}, but to the vector product of this equation with $\bx$, and also for  \eqref{Laws2phia}.\\
\fbox{}

In addition to the above equivalence result, it is easy to see that each individual term in \eqref{Laws2m}--\eqref{Laws2a} can be approximated arbitrarily well by the corresponding term  in \eqref{Laws2phim}--\eqref{Laws2phia}.
Indeed, fix any $\omega\subset\Omega$ with smooth $\dom$ and for sufficiently small $\eps>0$ define
\begin{equation} \label{eq:phieps}
	\phi_\eps=
	\begin{cases}
		\eps^{-1}\mbox{dist}(\bx,\dom),& \bx\in\mathcal{O}_\eps(\dom)\cap \omega,\\
		1& \bx\in\omega\setminus\mathcal{O}_\eps(\dom),\\
		0& \Omega\setminus \omega.
	\end{cases} 
\end{equation}
	We have $\phi_\eps\in W^{1,\infty}(\Omega)$ and one easily checks, letting $\tbn=-\nabla\phi_\eps/|\nabla\phi_\eps|$,  that
\begin{equation}\label{limits1}
	\frac{d}{dt}\int_{\omega}\phi_\eps\bu\,dx\to \frac{d}{dt}\int_{\omega}\bu\,ds, \quad 
	\int_{\omega}\bD(\bu)\tbn|\nabla\phi_\eps|\,dx\to \int_{\dom}\bD(\bu)\bn\,ds,~ \text{for}~\eps\to0,
\end{equation}
and smooth  $\bu$. Similarly, the limit values of other terms in  \eqref{Laws2phim}--\eqref{Laws2phia} will be their counterparts in \eqref{Laws2m}--\eqref{Laws2a}. Therefore, \eqref{Laws2phim}--\eqref{Laws2phia} can be also interpreted as the \emph{diffuse-volume} version of
conservation laws. Eqs. \eqref{Laws2phim}--\eqref{Laws2phia} imply   \eqref{Laws2m}--\eqref{Laws2a} term by term without the $\Div\bu=0$ assumption or equations \eqref{NSE} being invoked.

\medskip

Local conservation laws  in the Lagrangian form are written for the time-dependent material volume $\Omega_t$. Denote by $\Omega_0$ the fluid volume at a given initial moment $t=t_0$ and assume $\Omega_t\subset\Omega$ for $t\in [t_0,t_1]$ for some $t_1>t_0$. The evolution of $\Omega_t$ is defined by the Lagrangian mapping $\Phi_t:\Omega_0\to\Omega_t$, i.e. $\by=\Phi_t(\bx)$ solves the Cauchy problem 
\begin{equation}\label{Cauchy}
\by_t=\bu(t,\by),\quad t\in (t_0,t_1],\quad \by(t_0)=\bx.
\end{equation}

To properly reflect this domain evolution in
a weak form of \eqref{Laws1m}--\eqref{Laws1a}, we want $\phi$ to be time dependent and such that $\mbox{supp}(\phi)=\Omega_t$.
To this end, consider a smooth function $\phi^0$ such that $\mbox{supp}(\phi^0)=\Omega_0$. We define $\phi=\phi^0\circ\Phi_t^{-1}.$
The constructed $\phi$ is smooth (since $\bu$ is smooth  so is the solution to the Cauchy problem \eqref{Cauchy}), $\mbox{supp}(\phi)=\Omega_t$, and it satisfies the transport equation
\begin{equation}\label{eq:ls}
	\frac{\partial \phi}{\partial t}+(\bu\cdot \nabla)\phi
	=0\quad\text{in}~\Omega,~ t\in (t_0,t_1],\quad \phi(t_0)=\phi^0.
\end{equation}
Applying the Reynolds' transport theorem and using \eqref{eq:ls} and \eqref{NSE} one computes the following weak Lagrangian form of the local balances:
	\begin{align}
		\text{Moment.}\hskip7ex \frac{d}{dt}\int_{\Omega_t}\phi \bu\,dx&~=~ 2\nu\int_{\Omega_t}\bD(\bu)\tbn|\nabla\phi|\,dx-\int_{\Omega_t}p\tbn|\nabla\phi|\,dx, \label{Laws1phim}
		\\
		\text{Angl. Moment.}\quad  \frac{d}{dt}\int_{\Omega_t}\phi\bu\times\bx\,dx&
		~=~2\nu\int_{\Omega_t}(\bD(\bu)\tbn)\times\bx|\nabla\phi|\,dx-\int_{\Omega_t}p(\tbn\times\bx)|\nabla\phi|\,dx. \label{Laws1phia}
	\end{align}

By the same arguments as we use to prove Proposition~\ref{pr1} we prove the following proposition.
\begin{proposition}\label{pr2}
	Assume $\bu$ and $p$ are smooth and $\Div\bu=0$.  Then \eqref{Laws1m} {\rm (or \eqref{Laws1a})} holds for any material volume $\Omega_t$ such that $\Omega_t\subset\Omega$ for $t\in [t_0,t_1]$ iff  \eqref{Laws1phim} {\rm (or \eqref{Laws1phia})} holds for any $\phi$ satisfying \eqref{eq:ls} with  $\phi^0\in W^{1,\infty}(\Omega_{t_0})$, such that  $\mbox{supp}(\phi^0)=\Omega_{t_0}$.
\end{proposition}

Similar to the Eulerian case, it is easy to see that each individual term in \eqref{Laws1m}--\eqref{Laws1a} can be approximated arbitrarily well by the corresponding term  in \eqref{Laws1phim}--\eqref{Laws1phia}.
This time $\phi_\eps$ is constructed as $\phi_\eps=\phi^0_\eps\circ\Phi_t^{-1}$ with $\phi^0_\eps$ defined by the formula in \eqref{eq:phieps} with $\omega$ replaced by $\Omega_{t_0}$.
It holds $\phi_\eps\in W^{1,\infty}(\Omega\times[t_0,t_1])$ and one verifies, letting $\tbn=-\nabla\phi_\eps/|\nabla\phi_\eps|$,  that
\begin{equation}\label{limits2}
\frac{d}{dt}\int_{\Omega_t}\phi_\eps\bu\,dx\to \frac{d}{dt}\int_{\Omega_t}\bu\,ds, \quad 
\int_{\Omega_t}\bD(\bu)\tbn|\nabla\phi_\eps|\,dx\to \int_{\dO_t}\bD(\bu)\bn\,ds,~ \text{for}~\eps\to0,
\end{equation}
and smooth $\bu$. The limit values of other terms in  \eqref{Laws1phim}--\eqref{Laws1phia} will be their counterparts in \eqref{Laws1m}--\eqref{Laws1a}. Therefore, \eqref{Laws1phim}--\eqref{Laws1phia} can be also interpreted as the {diffuse-volume} version of local
conservation laws in the Lagrangian form.

\smallskip
In summary, equations \eqref{Laws2phim}--\eqref{Laws2phia} are equivalent formulations of the fundamental (local) conservation laws in the Eulerian formulation, while  \eqref{Laws1phim}--\eqref{Laws1phia} are equivalent formulations of the fundamental (local) conservation laws in the Lagrangian formulation. We will study the ability of a discretization method to match  \eqref{Laws2phim}--\eqref{Laws2phia} and \eqref{Laws1phim}--\eqref{Laws1phia}  instead of \eqref{Laws2m}--\eqref{Laws2a}  and  \eqref{Laws1m}--\eqref{Laws1a}.

\section{EMAC Galerkin formulation is locally conservative} \label{S4}
As an example of a continuous Galerkin method, we consider a conforming finite element method: Denote by $\bV_h\subset H^1_0(\Omega)^d$ and
$Q_h\subset L^2_0(\Omega)$ velocity and pressure   finite element spaces with respect to a tessellation $\T_h$ of $\Omega$ into elements (simplexes or more general polygons or polyhedra). We also need the following auxiliary spaces of continuous finite elements of degree $m+1$ and $m$, with $m\ge1$:
\begin{equation}\label{wbV}
	\begin{split}
		V_h&= \{v\in  H^1_0(\Omega):\,v\in\mathbb{P}_{m+1}(T)~\forall\,T\in\T_h\},\\ 
		\wbV&= \{v\in  H^1_0(\Omega):\,v\in\mathbb{P}_{m}(T)~~~~\forall\,T\in\T_h\}.
	\end{split}
\end{equation}

We only assume that the velocity space  contains all piecewise polynomial continuous functions of degree $m+1$, i.e. 
\begin{equation}\label{A1}
(V_h)^d\subset \bV_h.
\end{equation}
We do not have any further assumptions on finite element spaces, and in particular, both LBB stable and stabilized finite elements are admitted. 

\begin{remark}
\label{rem1} \rm Let  $\T_h$ be a triangulation  of $\Omega$ and $m\ge1$ be a polynomial degree.
The following examples of LBB stable FE pairs satisfy the assumption:  generalized Taylor--Hood $P_{m+1}-P_m$, $P_{m+1}-P_{m-1}^{\rm disc}$ (for $d=2$),  
$P_{m+1}-P_{m-2}^{\rm disc}$ (for $d=3$, $m>1$), $P_{m+1}^{\rm bubble}-P_{m-1}^{\rm disc}$ (for $d=3$ with face bubbles), generalized conforming Crouzeix--Raviart $P_{m+1}^{\rm bubble}-P_m^{\rm disc}$,
Scott-Vogelius $P_{m+1}-P_m^{\rm disc}$ (SV element is LBB stable subject to further assumptions on $\T_h$ \cite{guzman2018inf}), as well as LBB unstable  equal order  $P_{m+1}-P_{m+1}$ elements.
\end{remark}

We use $(f,g):=\int_\Omega f\cdot g\,dx$ notation for both scalar and vector functions $f,g$. The EMAC Galerkin formulation of \eqref{NSE} with $\bu={\bf 0}$ on $\dO$ reads: Find $\bu_h:(0,T)\to\bV_h$ and $\widehat p_h:(0,T)\to Q_h\cap L^2_0(\Omega)$ 
\begin{equation}\label{EMAC}
	\left\{
	\begin{aligned}
		\Big(\frac{\partial \bu_h}{\partial t},\bv_h\Big)+2(\bD(\bu_h)\bu_h,\bv_h)+ ((\Div\bu_h)\bu_h,\bv_h) +2\nu(\bD(\bu_h),\bD(\bv_h))+(\widehat p_h,\Div\bv_h) &=0\quad\forall\,\bv_h\in\bV_h, \\
		(\Div \bu_h,q_h)&=0 \quad\forall\,q_h\in Q_h,
	\end{aligned}\right.
\end{equation}
where $\widehat{p}_h$ approximates the EMAC pressure $\widehat{p}=p-\tfrac12|\bu|^2$. 
The EMAC formulation  is equivalent to other commonly used discrete formulations if  $\Div \bu_h=0$ pointwise.
However,  $(\Div \bu_h,q_h)=0$ does not imply $\Div \bu_h=0$ { except in special settings}.
 As a consequence, the discrete solution depends on the form of nonlinear terms used (i.e. EMAC, SKEW, CONV, ROT, etc.).
 
 Next, we demonstrate that the solution of \eqref{EMAC} satisfies discrete counterparts   of local conservation laws in both  Eulerian and Lagrangian forms. 
 
\subsection{Local conservation in Eulerian form} Unlike for the continuous problem, for the discrete case the counterparts of conservation laws in  Eulerian and Lagrangian forms do not follow one from another and we have to consider them separately. We start with the Eulerian form.
\medskip
  
\underline{Conservation of local linear momentum}. Consider arbitrary  $\phi_h\in V_h$, $\phi_h|_{\dO}=0$.  Then $\phi_h\be_i\in\bV_h$, for $i=1,\dots,d$, is a legitimate test function in \eqref{EMAC}. Letting $\bv_h=\phi_h\be_i$ in \eqref{EMAC} we compute for the nonlinear term
\begin{multline}
2(\bD(\bu_h)\bu_h,\phi_h\be_i)= (\bu_h\cdot\nabla\bu_h,\phi_h\be_i)+((\phi_h\be_i)\cdot\nabla\bu_h,\bu_h)\\
=
-(\bu_h\cdot\nabla(\phi_h\be_i),\bu_h)-((\Div\bu_h)\bu_h,\phi_h\be_i) -\tfrac12(\Div(\phi_h\be_i)\bu_h,\bu_h)
\\ =
-(\bu_h \cdot\nabla\phi_h,\bu_h\cdot\be_i)-((\Div\bu_h)\bu_h,\phi_h\be_i) - \tfrac12(\be_i\cdot\nabla\phi_h,|\bu_h|^2).
\end{multline}
Substituting this in the first equation from \eqref{EMAC} with  $\bv_h=\phi_h\be_i$ we obtain
\begin{equation*}
\Big(\frac{\partial \bu_h}{\partial t},\phi_h\be_i\Big)-(\bu_h \cdot\nabla\phi_h,\bu_h\cdot\be_i)- \tfrac12(\be_i\cdot\nabla\phi_h,|\bu_h|^2)
 +2\nu(\bD(\bu_h),\bD(\phi_h\be_i))-(\widehat p_h,\Div(\phi_h\be_i)) =0,
\end{equation*}
and after simple re-arrangements,
\begin{equation}\label{aux279}
\frac{d}{dt}\Big(\bu_h\cdot\be_i,\phi_h\Big)-(\bu_h \cdot\nabla\phi_h,\bu_h\cdot\be_i)+2\nu(\bD(\bu_h)\nabla\phi_h,\be_i)-(\widehat p_h+ \tfrac12|\bu_h|^2,\be_i\cdot\nabla\phi_h) =0.
\end{equation}
Let $\bn_h:=-\nabla\phi_h/|\nabla\phi_h|$ for $|\nabla\phi_h|\neq0$ (and arbitrary unit vector otherwise) and define
\[
\omega_h=\mbox{supp}(\phi_h)\quad \text{and}~~ p_h=\widehat p_h+ \tfrac12|\bu_h|^2,
\]
 then from equation \eqref{aux279}  for $i=1,\dots,d$ we get,
\begin{equation}\label{EmacLoc}
\frac{d}{dt}\int_{\omega_h}\phi_h\bu_h\,dx= 2\nu\int_{\omega_h}\bD(\bu_h)\bn_h|\nabla\phi_h|\,dx-\int_{\omega_h}p_h\bn_h|\nabla\phi_h|\,dx-\int_{\omega_h}\bu_h(\bu_h\cdot\bn_h)|\nabla\phi_h|\,dx,
\end{equation}
for any $\phi_h\in V_h$.
This is the discrete analogue of the local momentum conservation in \eqref{Laws2phim}.

\medskip
\underline{Conservation of local angular momentum}.  Consider arbitrary $\phi_h\in \wbV$.  Then  $\bx\times\phi_h\be_i\in\bV_h$ for $i=1,\dots,d$ is a legitimate test function. Letting $\bv_h=\bx\times\phi_h\be_i$ in \eqref{EMAC} we compute for the nonlinear term
\begin{multline}
	2(\bD(\bu_h)\bu_h,\bx\times\phi_h\be_i)= (\bu_h\cdot\nabla\bu_h,\bx\times\phi_h\be_i)+((\bx\times\phi_h\be_i)\cdot\nabla\bu_h,\bu_h)\\
	=
	-(\bu_h\cdot\nabla(\bx\times\phi_h\be_i),\bu_h)-((\Div\bu_h)\bu_h,\bx\times\phi_h\be_i) -\tfrac12(\Div(\bx\times\phi_h\be_i)\bu_h,\bu_h)
	\\ =
	-(\bu_h \cdot\nabla\phi_h,(\bu_h\times\bx)\cdot\be_i)-((\Div\bu_h)\bu_h,\bx\times\phi_h\be_i) - \tfrac12(\bx\times\nabla\phi_h,\be_i|\bu_h|^2).
\end{multline}
Substituting this in the first equation from \eqref{EMAC} with  $\bv_h=\bx\times\phi_h\be_i$ we obtain
\begin{multline}\label{aux324}
\Big(\frac{\partial \bu_h}{\partial t},\bx\times\phi_h\be_i\Big)	-(\bu_h \cdot\nabla\phi_h,(\bu_h\times\bx)\cdot\be_i)- \tfrac12(\bx\times\nabla\phi_h,\be_i|\bu_h|^2) \\
+2\nu(\bD(\bu_h),\bD(\bx\times\phi_h\be_i))-(\widehat p_h,\Div(\bx\times\phi_h\be_i)) =0.
\end{multline}
Simple re-arrangements give
\[
\frac{d}{dt}\Big(\bu_h\times\bx,\phi_h\be_i\Big)-(\bu_h \cdot\nabla\phi_h,(\bu_h\times\bx)\cdot\be_i) +2\nu(\bD(\bu_h)\nabla\phi_h,\bx\times\be_i)-((\widehat p_h+ \tfrac12|\bu_h|^2)\be_i,\nabla\phi_h\times\bx) =0.
\]
From the above equality for  $i=1,\dots,d$ we get
\begin{equation}\label{EmacLoc2}
	\frac{d}{dt}\int_{\omega_h}\phi_h\bu_h\times\bx\,dx =2\nu\int_{\omega_h}(\bD(\bu_h)\bn_h)\times\bx\,|\nabla\phi_h|\,dx-\int_{\omega_h}p_h(\bn_h\times\bx)|\nabla\phi_h|\,dx-\int_{\omega_h}(\bu_h\times\bx)(\bu_h\cdot\bn_h)|\nabla\phi_h|\,dx
\end{equation}
for any $\phi_h\in \wbV$.
This is the discrete analogue of the local angular momentum conservation from \eqref{Laws2phia}.

\begin{remark}\rm  \rev{ Conservation laws \eqref{EmacLoc} and \eqref{EmacLoc2} are local or element-wise in the sense that $\omega_h$ can be as small as the support of a nodal basis function from $V_h$ or $\wbV$, respectively. At the finite element level, they are no longer equivalent to standard element-wise conservation laws, such as the balances \eqref{Laws2m}--\eqref{Laws2a}, where $\omega$ and $\bu$ are replaced by $\omega_h$ and $\bu_h$. In particular, the argument in \eqref{limits1} and \eqref{limits2} is not valid at the finite element level; one cannot push $\epsilon$ to be smaller than $h$, which suggests an $O(h)$ discrepancy between the two formulations.
		
If convergence of $\bu_h$ and $p_h$ to the true smooth solution $\bu$ and $p$ is known, then an estimate of how accurate the finite element counterparts of \eqref{Laws2m}--\eqref{Laws2a} can be obtained through it. We are not pursuing such an estimate in this paper. {\color{blue}Instead}, the goal here is to formulate \textit{a priori} conservation laws for $\bu_h$ and $p_h$. This goal can be fulfilled {\color{blue} by} employing the weak forms of the conservation laws. The same remark remains largely valid for the element-wise balances in the Lagrangian form. 	 
}
\end{remark}


\subsection{Local conservation in Lagrangian form}
After discretization, there is no obvious equivalence between the Eulerian and Lagrangian forms of the local balances.  Nevertheless, one can show that EMAC form also obeys a discrete counterpart of the linear momentum local conservation in the Lagrangian form.  However, we need additional assumption on $V_h$ space.
Namely, we assume that the velocity space consists of piecewise polynomial continuous functions of degree $m+1$:
\begin{equation}\label{A2}
	(V_h)^d =\bV_h.
\end{equation}

\medskip
\underline{Conservation of local linear momentum}.
Consider $\phi_h^0\in V_h$ and 
$\phi_h:\,[t_0,\hat t_1]\to V_h$ solving
\begin{equation}\label{trans1}
\Big(\frac{\partial \phi_h}{\partial t},v_h\Big)+(\bu_h\cdot \nabla\phi_h,v_h)=0\quad\forall~v_h\in V_h,
\end{equation}
which is the projection of the transport equation \eqref{eq:ls} on the finite dimensional space $V_h$ with $\bu$ replaced by $\bu_h$. 

 Letting $\bv_h=\phi_h\be_i$ in \eqref{EMAC}, we repeat the same calculations as for the Eulerian case and arrive at \eqref{aux279}.
Since $\phi_h$ is time dependent, after re-arrangements  \eqref{aux279} gives  
\begin{equation*}
	\frac{d}{dt}\Big(\bu_h\cdot\be_i,\phi_h\Big)-\Big(\bu_h\cdot\be_i,\frac{\partial \phi_h}{\partial t}\Big)-(\bu_h \cdot\nabla\phi_h,\bu_h\cdot\be_i)+2\nu(\bD(\bu_h)\nabla\phi_h,\be_i)-(\widehat p_h+ \tfrac12|\bu_h|^2,\be_i\cdot\nabla\phi_h) =0.
\end{equation*}
Thanks to the assumption \eqref{A2} and equation \eqref{trans1}, the second and third terms add to zero.

Let $\bn_h:=-\nabla\phi_h/|\nabla\phi_h|$ for $|\nabla\phi_h|\neq0$ (and arbitrary unit vector otherwise) and define
\[
\Omega_h(t)=\mbox{supp}(\phi_h)\quad \text{and}~~ p_h=\widehat p_h+ \tfrac12|\bu_h|^2,
\]
then from equation \eqref{aux279}  for $i=1,\dots,d$ we get,
\begin{equation}\label{EmacLocL}
	\frac{d}{dt}\int_{\Omega_h(t)}\phi_h\bu_h\,dx= 2\nu\int_{\Omega_h(t)}\bD(\bu_h)\bn_h|\nabla\phi_h|\,dx-\int_{\Omega_h(t)}p_h\bn_h|\nabla\phi_h|\,dx,
\end{equation}
for any $\phi_h\in V_h$.
This is the discrete analogue of the local momentum conservation in \eqref{Laws1phim}.

\medskip
\underline{Conservation of local angular momentum}.
Consider $\phi_h^0\in \wbV$ and 
$\phi_h:\,[t_0,\hat t_1]\to \wbV$ solving
\[
\Big(\frac{\partial \phi_h}{\partial t},v_h\Big)+(\bu_h\cdot \nabla\phi_h,v_h)=0\quad\forall~v_h\in \wbV.
\] 

Letting $\bv_h=\bx\times\phi_h\be_i$ in \eqref{EMAC} we repeat the same calculations as for the Eulerian case and arrive at \eqref{aux324}.
Since $\phi_h$ is time dependent,  after re-arrangements \eqref{aux324} gives  
\begin{multline*}
\frac{d}{dt}\Big(\bu_h\times\bx,\phi_h\be_i\Big)-\Big(\frac{\partial \phi_h}{\partial t},(\bu_h\times\bx)\cdot\be_i\Big)-(\bu_h \cdot\nabla\phi_h,(\bu_h\times\bx)\cdot\be_i) \\ +2\nu(\bD(\bu_h)\nabla\phi_h,\bx\times\be_i)-((\widehat p_h+ \tfrac12|\bu_h|^2)\be_i,\nabla\phi_h\times\bx) =0.
\end{multline*}
Denote by $I_m(\bu_h\times\bx)$ a piecewise polynomial of degree $m$  interpolating $\bu_h\times\bx$, i.e. $I_m(\bu_h\times\bx)\in \wbV^3$. Then  $I_m(\bu_h\times\bx)\cdot\be_i\in\wbV$ holds.
Therefore, 
\[
\Big(\frac{\partial \phi_h}{\partial t},(\bu_h\times\bx)\cdot\be_i\Big)+(\bu_h \cdot\nabla\phi_h,(\bu_h\times\bx)\cdot\be_i)
=\left(\frac{d\phi_h}{d t},(\bu_h\times\bx-I_m(\bu_h\times\bx))\cdot\be_i \right)=: R_i,
\]
where $\frac{d\phi_h}{d t}=\frac{\partial \phi_h}{\partial t}+\bu_h \cdot\nabla\phi_h$.
Let $\bn_h:=-\nabla\phi_h/|\nabla\phi_h|$ for $|\nabla\phi_h|\neq0$ (and arbitrary unit vector otherwise) and define
\[
\Omega_h(t)=\mbox{supp}(\phi_h)\quad \text{and}~~ p_h=\widehat p_h+ \tfrac12|\bu_h|^2,
\]
then from equation \eqref{aux279}  for $i=1,\dots,d$ we get,
\begin{equation}\label{EmacLoc2L}
	\frac{d}{dt}\int_{\Omega_h(t)}\phi_h\bu_h\times\bx\,dx =2\nu\int_{\Omega_h(t)}(\bD(\bu_h)\bn_h)\times\bx\,|\nabla\phi_h|\,dx-\int_{\Omega_h(t)}p_h(\bn_h\times\bx)|\nabla\phi_h|\,dx+R
\end{equation}
for any $\phi_h\in \wbV$.
This is the discrete analogue of the local momentum conservation in \eqref{Laws1phia} \emph{up to the residual term} $R=R_1+\dots+R_d$. If we assume that $\bu_h$ approximates a (smooth) solution to the NSE with order $O(h^{r})$, $r\ge m+1$, in some norm $\|\cdot\|_\ast$ then $\|R\|_\ast=O(h^{m+1})$ once $\frac{d\phi_h}{d t}$ is bounded in the dual norm to $\|\cdot\|_\ast$. According to \eqref{A2} the optimal approximation order for  $\bu_h$ would be $O(h^{m+2})$ in $L^2(L^2)$ norm.

\section{Numerical Tests} \label{S5}

We now give numerical examples to illustrate the theory above.  For these tests, the full Navier--Stokes discretization uses BDF temporal discretizations (and Crank--Nicolson for the initial time steps) and $(\bV_h,Q_h)$ \rev{ is the $P_2-P_1$} Taylor--Hood elements on a mesh \rev{$\T_h$}.  {\color{red}
The schemes used to compute solutions are (at time step $n$): Find $(\bu_h^n,P_h^n)\in (\bV_h,Q_h)$ satisfying
\begin{align*}
	\left( \left(\frac{d \bu_h}{dt}\right)^{n}_{\rm bdfk},\bv_h \right) + 2((\bD(\bu_h^n)\bu_h^n,\bv_h) - (p_h^n,\nabla \cdot \bv_h) + 2\nu(\bD(\bu_h^n),\bv_h)&= {\bf f(t^n)}, \\
(\nabla \cdot \bu_h^n,q_h)& =0,
\end{align*}
for all $(\bv_h,q_h)\in (\bV_h,Q_h)$.  The BDF notation for the time derivative term used above is defined as follows.  For a sequence $\{f^n\}_{n=0,1,\dots}$ of scalar or vector quantities (where $n$ denotes a time level), and $\Delta t$ is denoting a time step size, we use the shortcut notations for discrete time derivatives:
	\begin{align*}
	\left(\frac{df}{dt}\right)^{n}_{\rm bdf3}& = \frac{\frac{11}{6} f^{n} - 3 f^{n-1} + \frac{3}{2} f^{n-2} - \frac13 f^{n-3} }{\Delta t}, \\
		\left(\frac{df}{dt}\right)^{n}_{\rm bdf2}&= \frac{3 f^{n} - 4f^{n-1} + f^{n-2}}{2\Delta t}, \\
		 \left(\frac{df}{dt}\right)^{n}_{\rm bdf1}&= \frac{f^{n} - f^{n-1}}{\Delta t}.
	\end{align*}
For our computations below, k=2 or 3 for the NSE schemes and k=1 or 2 for the discrete transport equations.}. The nonlinear problem at each time step is resolved with Newton's method, and typically it takes just 2 or 3 iterations to resolve.  

With temporal discretizations, the precise definitions of the discrete local balances will change accordingly, and we derive these now before proceeding to the tests.  Denote by $\omega_h$ the approximation of a subdomain $\omega$ whose boundary consists of element edges from the mesh.  Define the functions $\phi_h\in V_h$ and $\psi_h \in \tilde V_h$ nodally by
\begin{align}
\phi_h(x_j) = \left\{\begin{array}{c} 1 \ \mbox{ if $x_j$ is a node on $P_2(\tau_h)$ in the interior of $\omega_h$} \\0 \ \mbox{ otherwise,         } \ \ \ \ \ \ \ \ \ \ \ \ \ \ \ \ \ \ \ \ \ \ \ \ \ \ \ \ \ \ \ \ \ \ \ \ \ \ \  \ \ \ \ \   \end{array}\right. \label{phi1} \\
\psi_h(x_j) = \left\{\begin{array}{c} 1 \ \mbox{ if $x_j$ is a node on $P_1(\tau_h)$ in the interior of $\omega_h$} \\0 \ \mbox{ otherwise.         } \ \ \ \ \ \ \ \ \ \ \ \ \ \ \ \ \ \ \ \ \ \ \ \ \ \ \ \ \ \ \ \ \ \ \ \ \ \ \ \ \ \ \ \    \end{array}\right. \label{psi1}
\end{align}
In our implementations we apply BDF formulas for the temporal discretization for the momentum and transport equations.

\subsubsection{Discrete local conservation in Eulerian form}

We consider first the discrete Eulerian form of local conservation of momentum and angular momentum.
Choosing $\phi_h$ by \eqref{phi1} and repeating the arguments above that derived \eqref{EmacLoc} but using the BDFk temporal discretization, we get the following (fully) discrete local momentum balance
\begin{align*}
\left(\frac{d\Big(\int_{\omega_h}  \phi_h\bu_h\,dx\Big)}{dt}\right)^{n}_{\rm bdfk} = 2\nu\int_{\omega_h}\bD(\bu_h^{n})\bn_h|\nabla\phi_h|\,dx-\int_{\omega_h}p_h^{n}\bn_h|\nabla\phi_h|\,dx-\int_{\omega_h}\bu_h^{n}(\bu_h^{n}\cdot\bn_h)|\nabla\phi_h|\,dx,
\end{align*}
with $\bn_h=-\nabla\phi_h/|\nabla\phi_h|$. 
Similarly, for discrete local angular momentum conservation we obtain
\begin{align*}
	\left(\frac{d\Big(\int_{\omega_h}  \psi_h\bu_h\times\bx\,dx\Big)}{dt}\right)^{n}_{\rm bdfk}
	 =&2\nu\int_{\omega_h}(\bD(\bu_h^{n})\bn_h)\times\bx\,|\nabla\psi_h|\,dx\\
	 &-\int_{\omega_h}p_h^{n}(\bn_h\times\bx)|\nabla\psi_h|\,dx-\int_{\omega_h}(\bu_h^{n}\times\bx)(\bu_h^{n}\cdot\bn_h)|\nabla\psi_h|\,dx,
\end{align*}
where $\psi_h$ is defined by \eqref{psi1}, and $\bn_h=-\nabla\psi_h/|\nabla\psi_h|$.

In our tests, we will show plots of discrete local Eulerian momentum error
\begin{multline*}
e_{E}^{mom} = 
\left(\frac{d\Big(\int_{\omega_h}  \phi_h\bu_h\,dx\Big)}{dt}\right)^{n}_{\rm bdfk}  -  2\nu\int_{\omega_h}\bD(\bu_h^{n})\bn_h|\nabla\phi_h|\,dx\\ +\int_{\omega_h}p_h^{n}\bn_h|\nabla\phi_h|\,dx+\int_{\omega_h}\bu_h^{n}(\bu_h^{n}\cdot\bn_h)|\nabla\phi_h|\,dx ,
\end{multline*}
and discrete local Eulerian angular momentum error
\begin{multline*}
e_{E}^{am} 	= 
\left(\frac{d\Big(\int_{\omega_h}  \psi_h\bu_h\times\bx\,dx\Big)}{dt}\right)^{n}_{\rm bdfk}   -  2\nu\int_{\omega_h}(\bD(\bu_h^{n})\bn_h)\times\bx\,|\nabla\psi_h|\,dx\\ +\int_{\omega_h}p_h^{n}(\bn_h\times\bx)|\nabla\psi_h|\,dx+\int_{\omega_h}(\bu_h^{n}\times\bx)(\bu_h^{n}\cdot\bn_h)|\nabla\psi_h|\,dx .
\end{multline*}

\subsubsection{Discrete local conservation in Lagrangian form}

Discrete local conservation in Lagrangian form is somewhat more complicated compared to the Eulerian case due to the $\phi_h$ function becoming time dependent in the momentum and angular momentum balances, as well as the transport equations involved in these balances being hyperbolic.  As our tests are for illustrative purposes of certain theoretical properties, we approximate the transport equations in the following way for the purpose of ease in computations, even though
other approaches to solving the transport equation may be better in practice.  

Consider $\phi_h^0\in V_h$ to be defined by \eqref{phi1}, and then define $\phi_h^{n} \in V_h$ (n=1,2,3,...) via
\begin{equation}\label{trans1h}
\Big(	\left(\frac{d\phi_h}{dt}\right)^{n}_{\rm bdfj}  
 ,v_h\Big)+(\bu_h^{n} \cdot \nabla\phi_h^{n} ,v_h) 
 =0\quad\forall~v_h\in V_h,
\end{equation}
where $j=$1 or 2 in our numerical tests (and if $j=$2 then  the first time step is backward Euler).


Rederiving the discrete local Lagrangian momentum balance \eqref{EmacLocL} but now using BDFk ($k=$2 or 3) time stepping for Navier-Stokes together with \eqref{trans1h}, we obtain the balance
\begin{equation*}
\int_{\Omega}  \left( \phi_h 	\left(\frac{d\bu_h}{dt}\right)^{n}_{\rm bdfk}  +
	\left(\frac{d\phi_h}{dt}\right)^{n}_{\rm bdfj}  \bu_h^{n}\right) \ dx 
 = 2\nu\int_{\Omega }\bD(\bu_h^{n})\bn^{n}_h|\nabla\phi^{n}_h|\,dx-\int_{\Omega}p_h^{n} \bn^{n}_h|\nabla\phi^{n}_h|\,dx,
\end{equation*} 
and thus define the discrete local Lagrangian momentum error by
\begin{equation*}
e_L^{mom} = 
\int_{\Omega}  \left( \phi_h^{n} 	\left(\frac{d\bu_h}{dt}\right)^{n}_{\rm bdfk} +
	\left(\frac{d\phi_h}{dt}\right)^{n}_{\rm bdfj}  \bu_h^{n}\right) \ dx 
 - \left( 2\nu\int_{\Omega }\bD(\bu_h^{n})\bn^{n}_h|\nabla\phi^{n}_h|\,dx-\int_{\Omega}p_h^{n} \bn^{n}_h|\nabla\phi^{n}_h|\,dx \right).
\end{equation*} 
Note that if the transport equations are solved in a different way, then these definitions of discrete Lagrangian momentum and angular momentum balances need modified accordingly.  For example if an explicit method is used, then the local balance will be defined with some terms at time $t^{n}$ and others and time $t^{n-1}$.

For angular momentum, we proceed similarly as for momentum to find a fully discrete analogue to \eqref{EmacLoc2L}.  Let $\psi_h^0$ be defined by \eqref{psi1} and find $\psi_h^{n}\in \tilde V_h\cap H^1_0(\Omega)$ for $n=1,2,3,...$ by 
\[
\Big(\left(\frac{d\psi_h}{dt}\right)^{n}_{\rm bdfj},v_h\Big)+(\bu_h^{n} \cdot \nabla\psi_h^{n},v_h) =0\quad\forall~v_h\in \tilde V_h\cap H^1_0(\Omega),
\] 
and using backward Euler for the first time step if $j=$2.  Following similar steps as the theory above, the discrete Lagrangian local angular momentum error is then 
given by 
\begin{align*}
 e_L^{am}  = \bigg( \int_{\Omega}  \bigg( \psi_h^{n} &   	\left(\frac{d\bu_h}{dt}\right)^{n}_{\rm bdfk}   + \left(\frac{d\psi_h}{dt}\right)^{n}_{\rm bdfj} \bu_h^{n}  \bigg) \times \bx 
   \, dx
	\\
	& -  2\nu\int_{\Omega}(\bD(\bu_h^{n})\bn_h^{n})\times\bx\,|\nabla\psi_h^{n}|\,dx+\int_{\Omega}p_h^{n}(\bn_h^{n}\times\bx)|\nabla\psi_h^{n}|\,dx \bigg).
\end{align*}

{\color{red}
\subsubsection{Discrete local conservation in traditional form}

To illustrate the non-equivalence at the discrete level of the traditional local conservation and the proposed weak formulations, we consider also the momentum and angular momentum error from using the discrete solutions $(\bu_h,p_h)$ in Eulerian conservation laws \eqref{Laws2m} and \eqref{Laws2a}, and approximating the time derivative of the velocity with the BDF approximation used in that simulation.  Hence we define errors in traditional discrete Eulerian local conservation by
\[
e^{mom}_{trad} =  \int_{\omega}\left( \frac{d \bu_h}{dt}\right)^n_{\rm bdfk} \,dx - \left( 2\nu\int_{\dom}\bD(\bu_h^n)\bn\,ds - \quad\int_{\dom}p_h^n\bn\,ds-\int_{\dom}\bu_h^n(\bu_h^n\cdot\bn)\,ds \right),
\]
and
\[
e^{am}_{trad} = 
\int_{\omega} \left( \frac{d \bu_h}{dt}\right)^n_{\rm bdfk} \times\bx\,dx
 -\left( 2\nu\int_{\dom}(\bD(\bu_h^n)\bn)\times\bx\,ds
 - \int_{\dom}p_h^n(\bn\times\bx)\,ds
 - \int_{\dom}(\bu_h^n \times\bx)(\bu_h^n \cdot\bn)\,ds \right).
 \]

We note we could also consider errors of conservation in discrete traditional Lagrangian form.  Doing this, however, is a difficult computational task in most finite element codes, and so we omit this comparison.

}

\subsection{Gresho problem}

For our first test we use a slight variation of the classical Gresho problem  on $\Omega=(-0.5,0.5)^2$, which consists 
\begin{wrapfigure}{r}{0.37\textwidth}
	\centering
		\vskip-2ex
	\includegraphics[width=0.35\textwidth,height=0.35\textwidth, trim=0 0 0 0, clip]{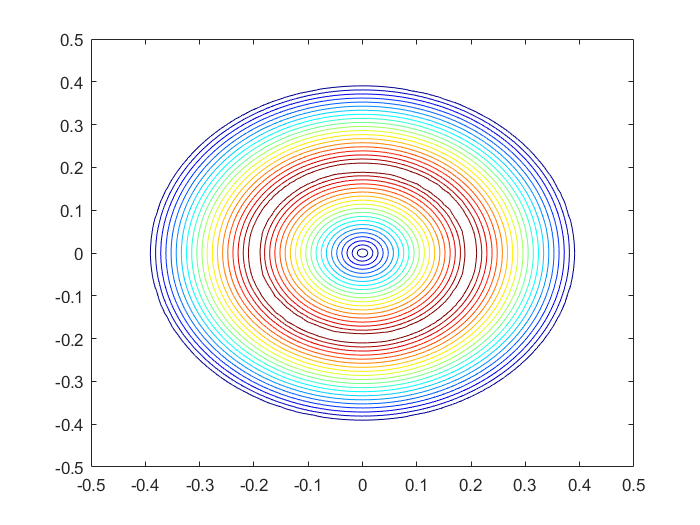}
	\vskip-3ex
	\caption{\small Initial velocity for the Gresho problem is shown above, as speed contours.}
	\label{fig:True_Gresho}
		\vskip-18ex
\end{wrapfigure}
of a velocity and pressure
\begin{align*}
	\small
\bu=\begin{cases}
\bmat{-5y~5x}^T &\text{for }r< 2,\\[1ex]
\bmat{\frac{2y}{r}+5y~\frac{2x}{r}-5x}^T &\text{for }.2 \leq r\leq .4,\\
\bmat{0~0}^T &\text{for }r>.4,
\end{cases},
\end{align*}
\begin{align*}
	\small
p=\begin{cases}
12.5r^2+C_1 &\text{for }r<.2,\\
12.5r^2-20r+4\log(r)+C_2 &\text{for }.2 \leq r \leq .4,\\
0 &\text{for }r>.4,
\end{cases}
\end{align*}
where $r=\sqrt{x^2+y^2}$ and 
{\small
\begin{align*}
&C_2=-12.5(.4)^2+20(.4)^2-4\log(.4),\\
&C_1=C_2-20(.2)+4\log(.2).
\end{align*}
}
This velocity is plotted in figure \ref{fig:True_Gresho} and is an exact solution of the unforced steady Euler equations, and hence an accurate solver should preserve the initial condition in time.  It is shown in \cite{CHOR17, OR20} that a NSE solver with EMAC nonlinearity and using Crank-Nicolson time stepping together with Taylor-Hood finite element spatial discretization will preserve pointwise global energy, momentum and angular momentum for this problem while other common nonlinearity formulations such as SKEW, ROT and CONV will not preserve these physical balance laws and moreover will be less accurate in the sense of $L^2(\Omega)$ error.

\begin{figure}[h]
	\centering
	\includegraphics[width=0.3\textwidth,height=0.28\textwidth, trim=0 0 0 0, clip]{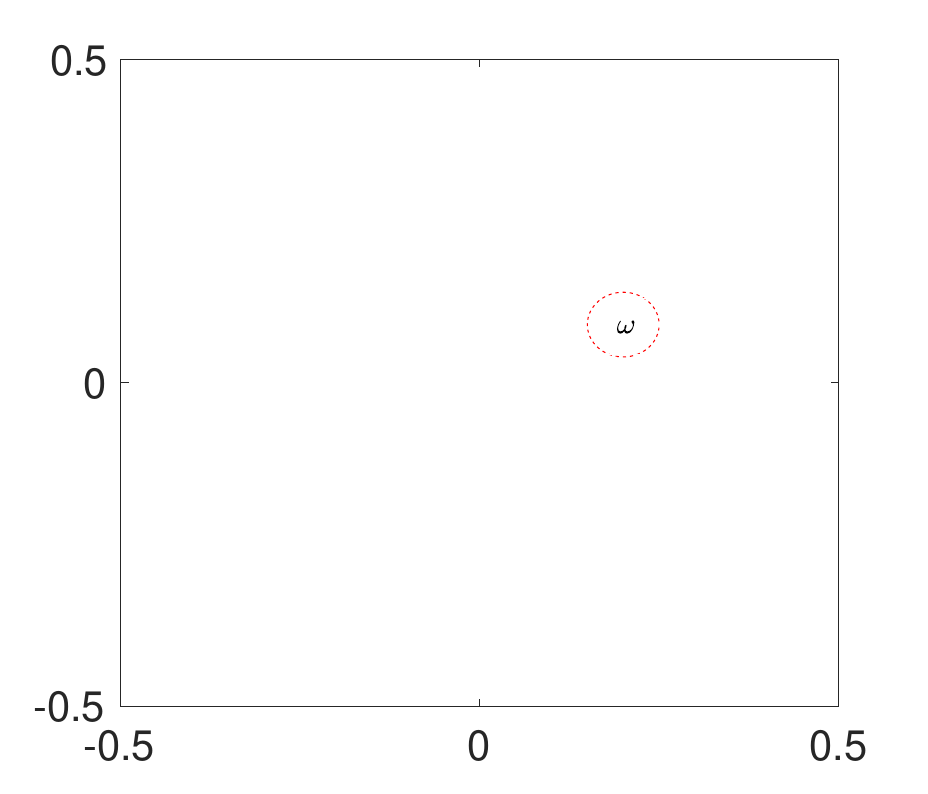}
	\includegraphics[width=0.3\textwidth,height=0.28\textwidth, trim=0 0 0 0, clip]{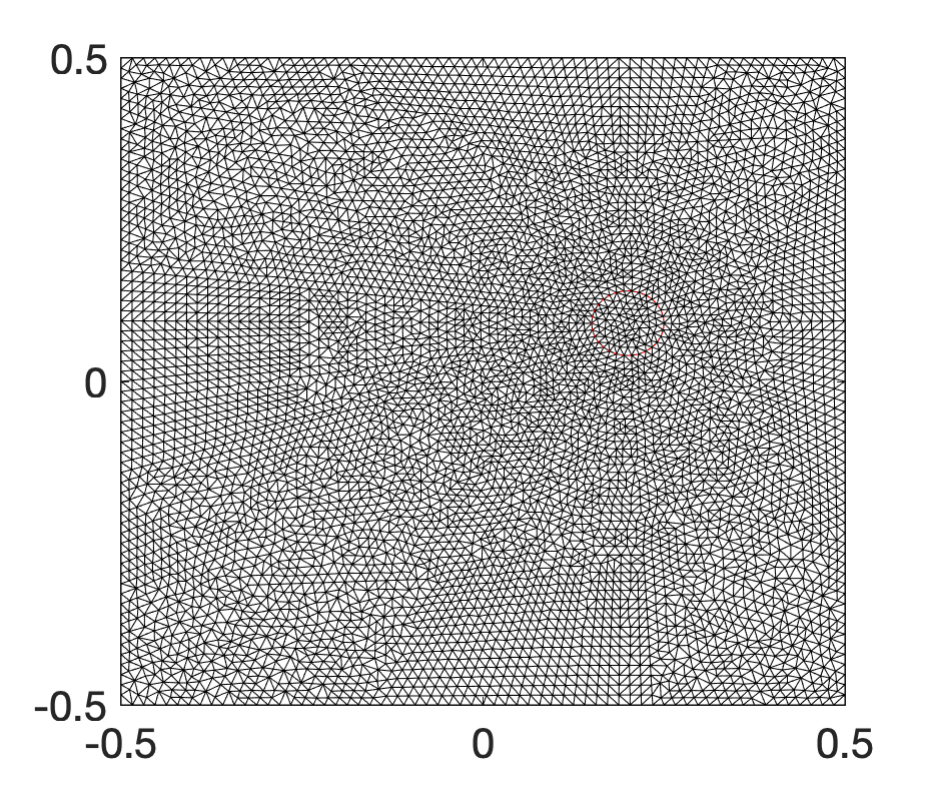}
	\includegraphics[width=0.3\textwidth,height=0.28\textwidth, trim=0 0 0 0, clip]{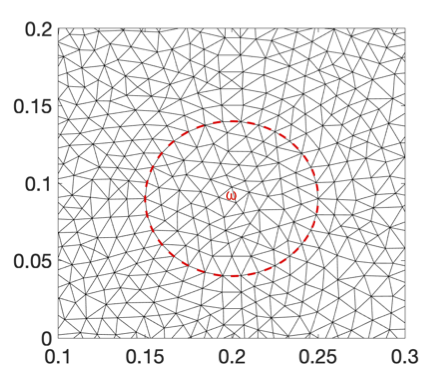}
	\caption{Shown above is the domain and $\omega$ (left), the mesh (center), and the mesh zoomed in near $\omega$ for the Gresho problem.}
	\label{Gresho}
\end{figure}

{\color{red} 
We alter this problem very slightly by changing the viscosity to $\nu=10^{-10}$ so as not to solve the Euler equations but instead the NSE.  We note this change of viscosity will (very slightly) change the true solution in time, however this is of no consequence as our interest herein is not the solution but the local conservation of momentum and angular momentum.  We choose $\omega$ to be the circle of radius $0.05$ centered at $(0.2,0.09)$, as shown in figure \ref{Gresho} at left.  Figure \ref{Gresho} also shows the mesh $\tau_h$ used for the computations below as well as the mesh zoomed in near $\omega$; the mesh is a Delaunay triangulation constructed from having 65 nodes on each domain edge and 30 nodes on $\partial\omega$.   We define $\omega_h$ to be the approximation of $\omega$ whose boundary consists of triangle edges from the mesh.

\begin{figure}[h]
	\begin{center}
		\includegraphics[width=0.32\textwidth,height=0.27\textwidth, trim=0 0 40 0, clip]{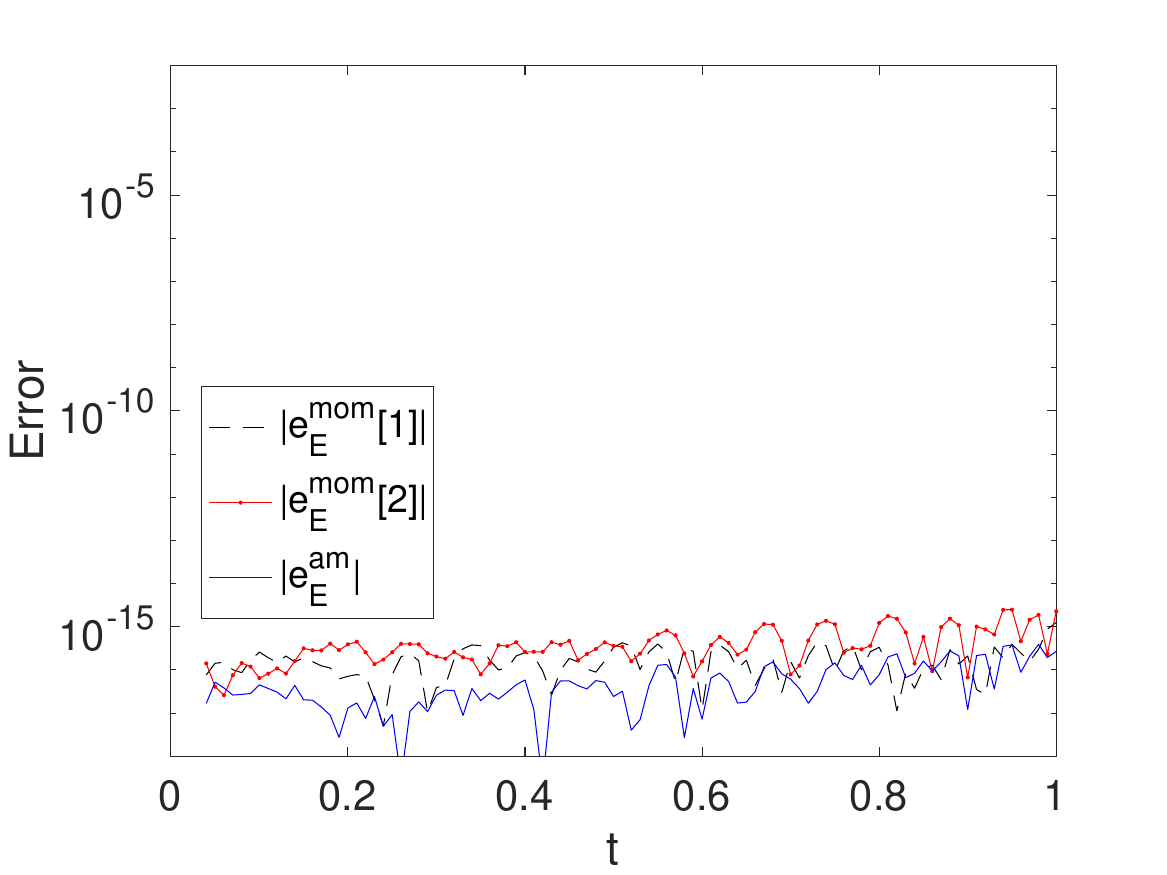}
		\includegraphics[width=0.32\textwidth,height=0.27\textwidth, trim=0 0 40 0, clip]{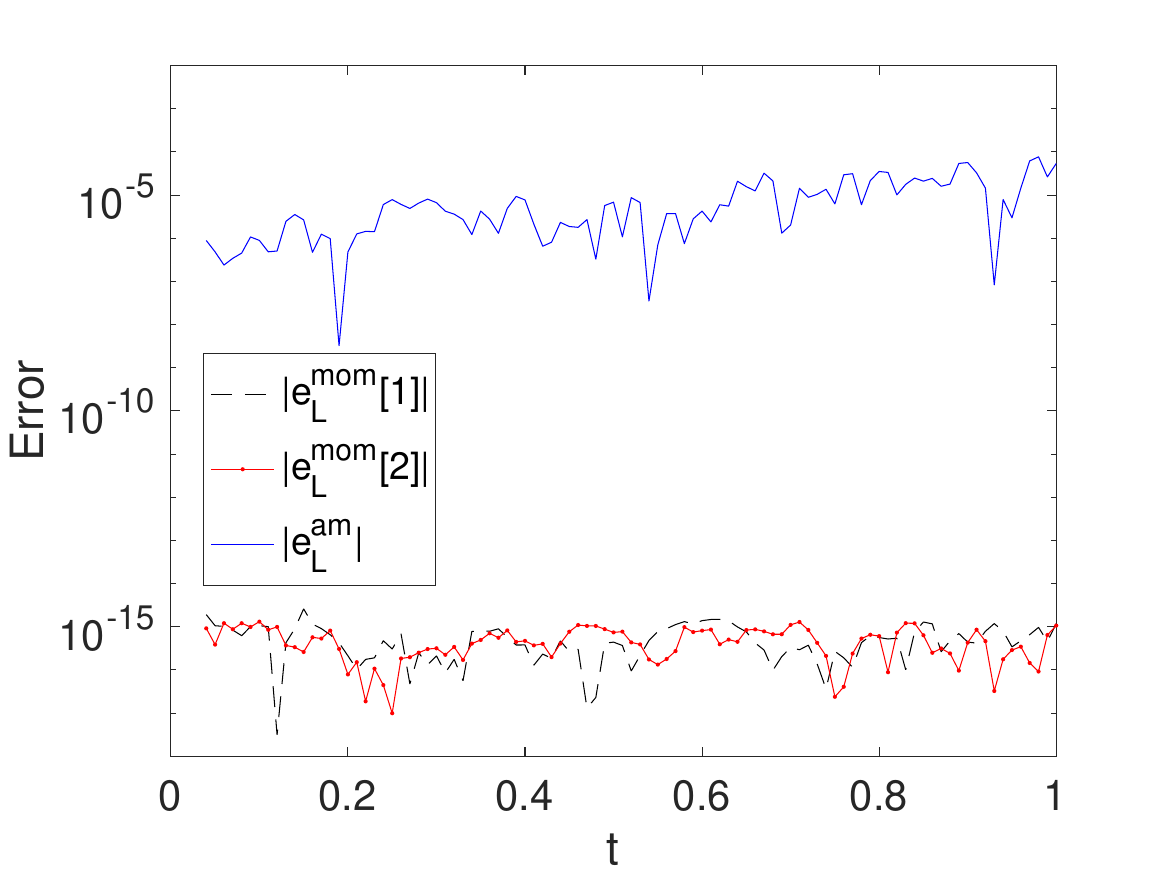}
		\includegraphics[width=0.32\textwidth,height=0.27\textwidth, trim=0 0 40 0, clip]{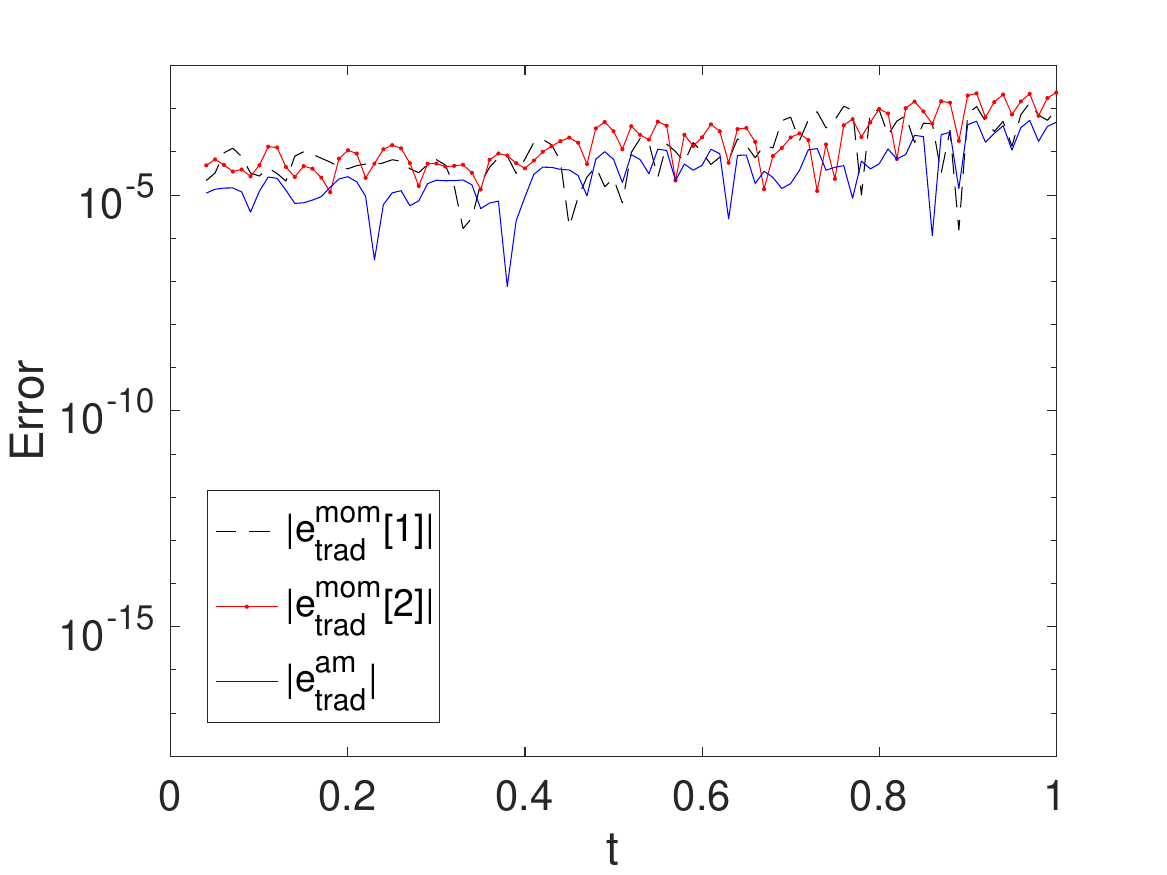}
	\end{center}
	\caption{Shown above is error {\color{red} in discrete local Eulerian (left), Lagrangian (center), and traditional Eulerian (right) momentum and angular momentum conservation }versus time in the (viscous) Gresho problem.  \label{EerrG} }
\end{figure}

Computations are done using end time $T=1$, time step size $\Delta t=0.01$, no external forcing, initial condition $\bu_h^0$ is the nodal interpolant of the true solution, and we show errors in discrete Eulerian, Lagrangian and traditional Eulerian conservation in figure \ref{EerrG} as absolute values of errors versus time. Notation $\mathrm{e}^{\rm mom}_{\rm E}[i]$ is used for the $i$-th component of the linear momentum.   We observe these quantities are conserved pointwise for Eulerian, and are stable in time, just as the theory above predicts.   Discrete local Lagrangian momentum is also preserved pointwise, although Lagrangian angular momentum is not preserved pointwise but instead has values as large as $O(10^{-4})$, which is consistent with the $O(h^2)$ residual our theory above predicts.  Discrete traditional Eulerian momentum and angular momentum are conserved up to about $O(10^{-3})$ or so, which is consistent with the spatial discretization error being $O(h^2)$ in the gradient of the velocity.}

\subsection{2D flow past a cylinder}

\begin{figure}[h!]
\centering
\includegraphics[width=0.98\textwidth,height=0.25\textwidth, trim=120 0 120 0, clip]{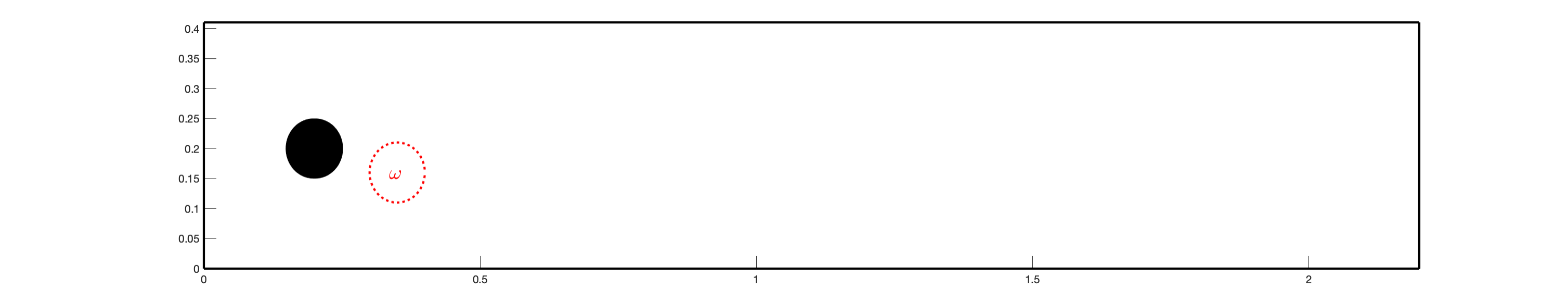} \\
\caption{The domain for the channel flow past a cylinder numerical experiment.\label{cyldomain} }
\end{figure}

\begin{figure}[h]
	\begin{center}
		\includegraphics[width=0.48\textwidth,height=0.15\textwidth, trim=80 10 60 10, clip]{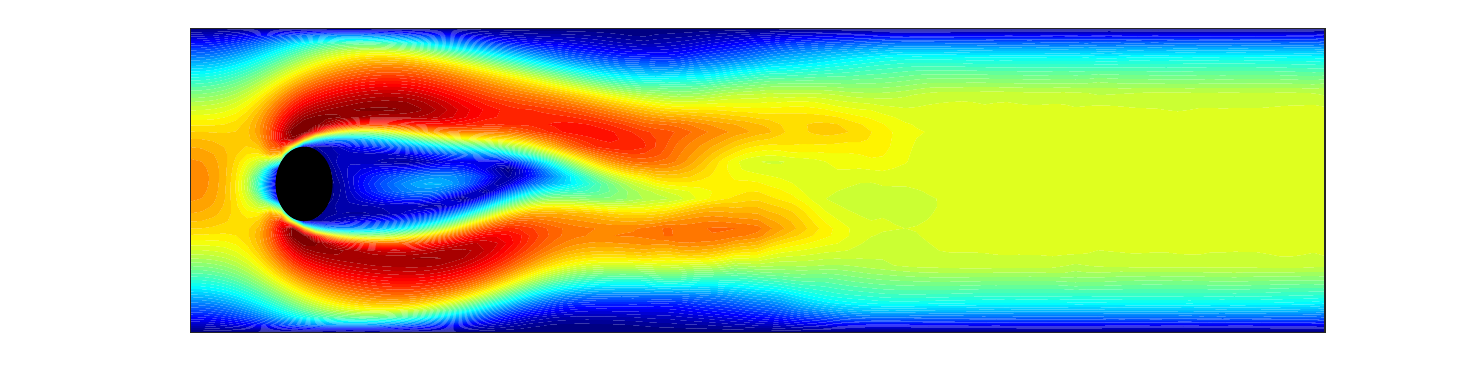}
		\includegraphics[width=0.48\textwidth,height=0.15\textwidth, trim=80 10 60 10, clip]{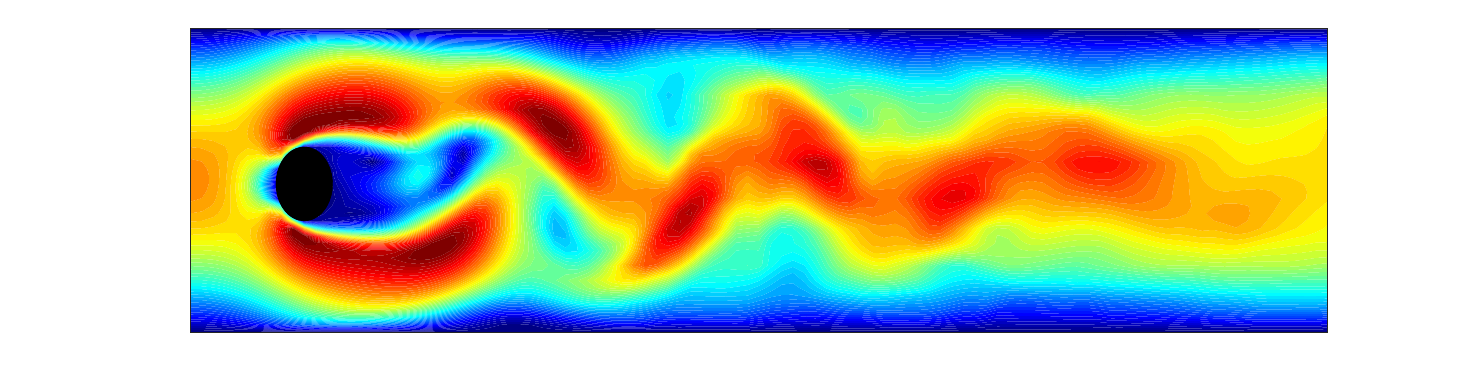}
		\includegraphics[width=0.48\textwidth,height=0.15\textwidth, trim=80 10 60 10, clip]{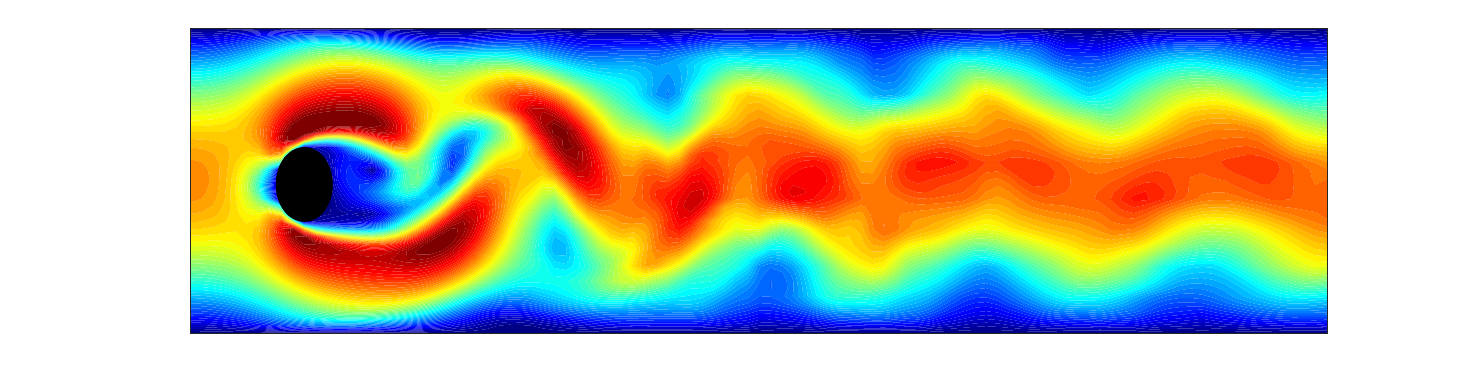}
		\includegraphics[width=0.48\textwidth,height=0.15\textwidth, trim=80 10 60 10, clip]{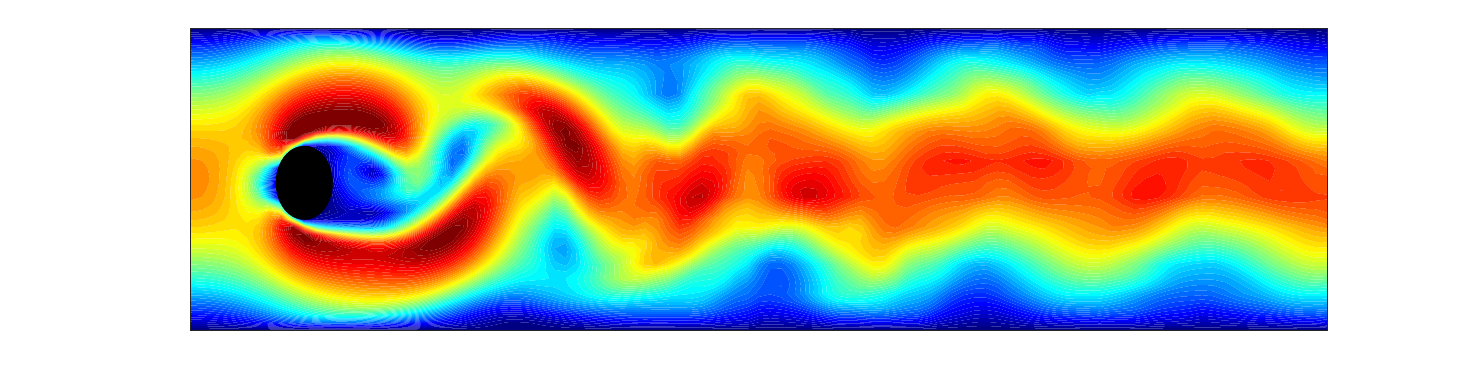}
		\includegraphics[width=0.48\textwidth,height=0.15\textwidth, trim=80 10 60 10, clip]{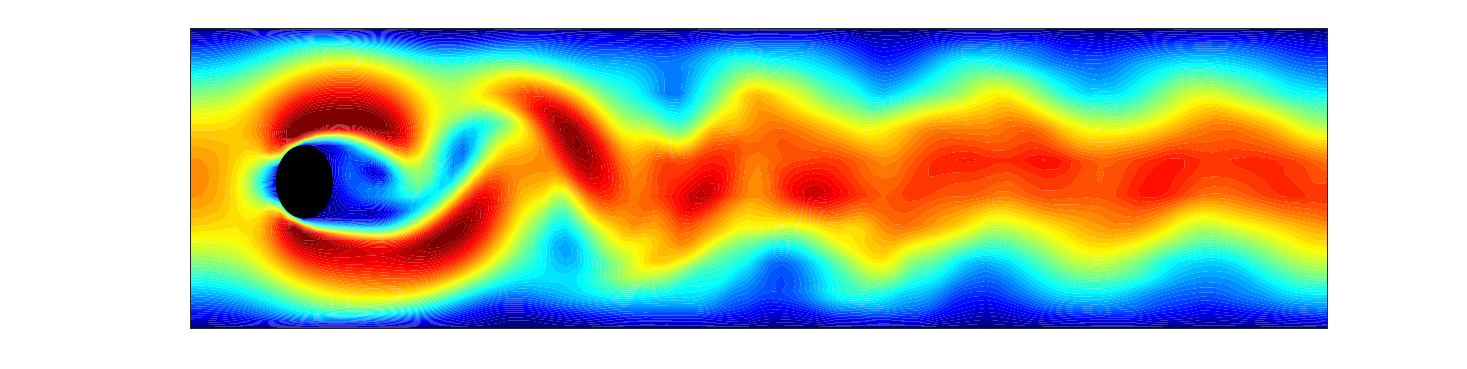}
	\end{center}
	\caption{Shown above are the t=1,2,3,4,5 solution plots of the Re=100 simulations of flow past a cylinder, as speed contours. \label{cylplots} }
\end{figure}

\begin{figure}[h!]
\begin{center}
\includegraphics[width=0.4\textwidth,height=0.25\textwidth, trim=0 0 0 0, clip]{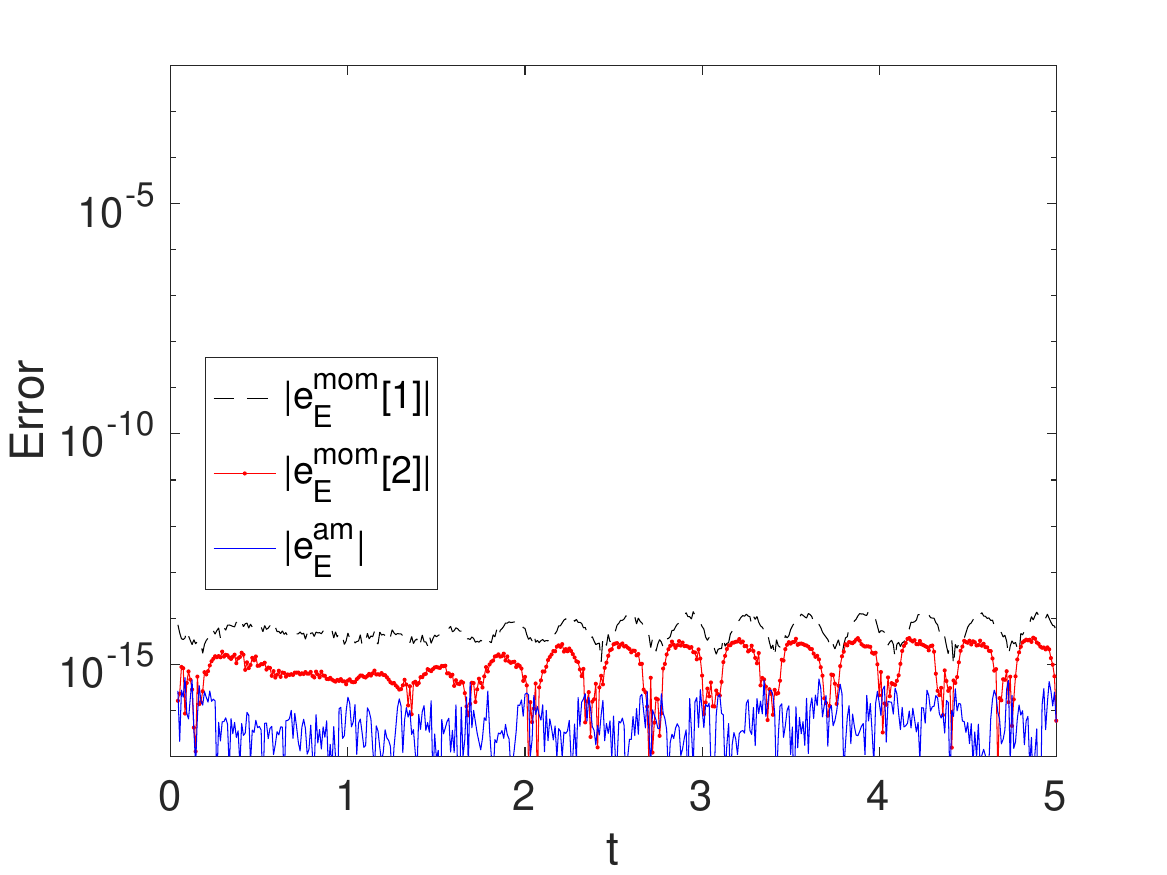}
\includegraphics[width=0.4\textwidth,height=0.25\textwidth, trim=0 0 0 0, clip]{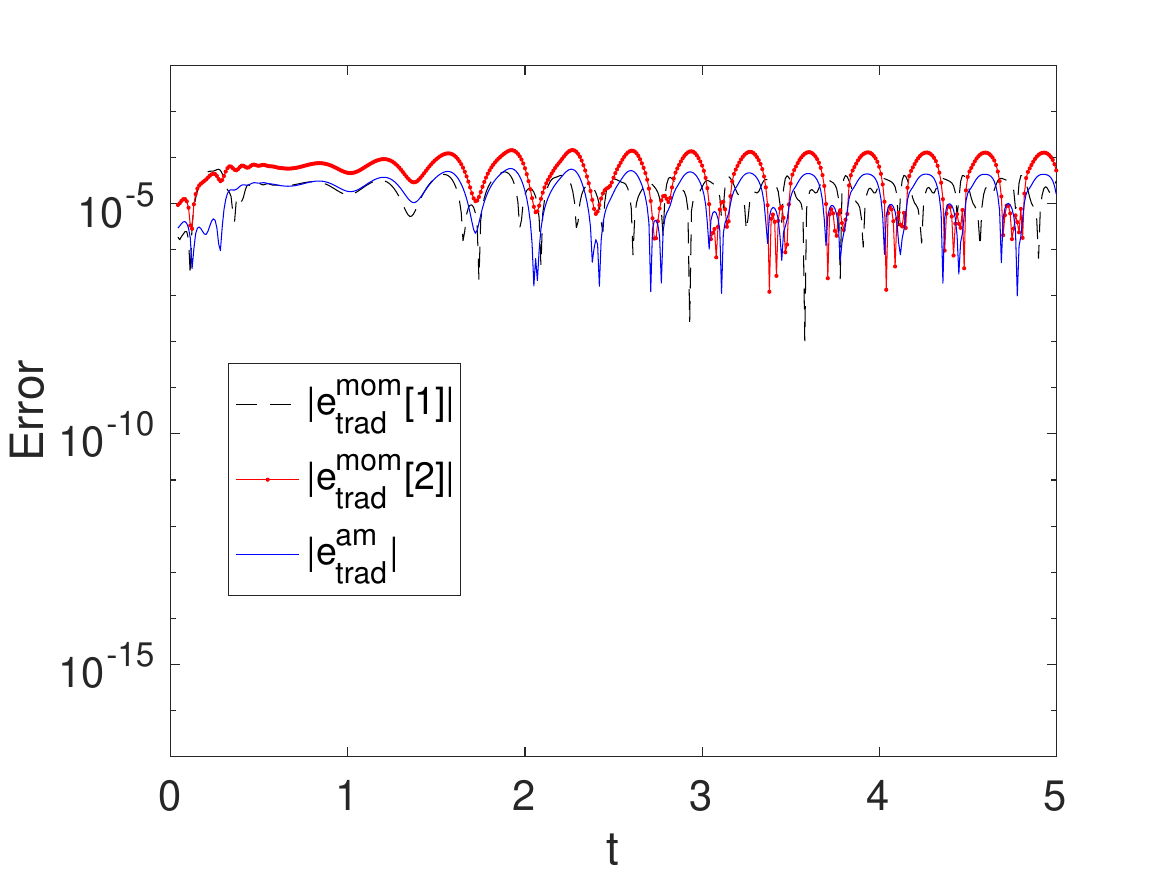}
\end{center}
\caption{Shown above is error in discrete local Eulerian {\color{red} (left) and traditional discrete local Eulerian (right)} momentum and angular momentum conservation versus time in the 2D channel flow past a cylinder test.\label{cylerr} }
\end{figure}

{\color{red}
For our next test we consider the classical 2D channel flow past a cylinder problem originally from \cite{ST96}, but with updated benchmark data and descriptions in \cite{J04,J16} and references therein.  The domain is the rectangle $[0,2.2]\times [0,0.41]$ as shown in figure \ref{cyldomain}, with a cylinder centered at $(0.2,0.2)$ with radius $0.05.$ We take no external forcing, set $\nu=0.001$ (which corresponds to Reynolds number 100, using the mean inflow velocity of 1), and set inflow/outflow profiles to be
\begin{align*}
u_1(0,y,t) & = u_1(2.2,y,t) = \frac{6}{0.41^2}y(0.41-y),
\\u_2(0,y,t) & = u_2(2.2,y,t) = 0.
\end{align*}
The flow starts from rest, and solution plots at times t=1,2,3,4,5 from our computations described below are shown in figure \ref{cylplots}, as speed contours.  By t=4, the flow has reached a periodic in time state with the repeating Van Karman vortex street.}

We define a subdomain $\omega$ to be a circle radius $0.05$ centered at $(0.35,0.16)$, and $\omega_h$ to be its approximation by the mesh.  A plot of $\omega$ is shown in figure \ref{cyldomain}.  We use a mesh that provides approximately 64K velocity degrees of freedom (dof) and 7K pressure dof when discretized with Taylor-Hood elements.  We compute using BDF3 time stepping {\color{red} to $T=5$ using time step size $\Delta t=0.01$, and start the flow from rest, $\bu_0 = {\bf 0}$.  Errors in discrete local Eulerian momentum and angular momentum conservation are shown in figure \ref{cylerr} and we once again observe pointwise local conservation.  Traditional Eulerian momentum and angular momentum are also shown, and we observe those to be $O(10^{-5})$, which is consistent with the discretization error.  We do not consider Lagrangian discrete local conservation for this test, since there is an outflow and conservation is therefore not expected except for very short times as the transported quantity will exit the domain through the outflow.  We note also that all results from this test are very similar if BDF2 is used instead of BDF3.}

\subsection{Kelvin-Helmholtz flow}

\begin{figure}[h]
	\begin{center}
		\includegraphics[width=0.3\textwidth,height=0.27\textwidth, trim=0 0 0 0, clip]{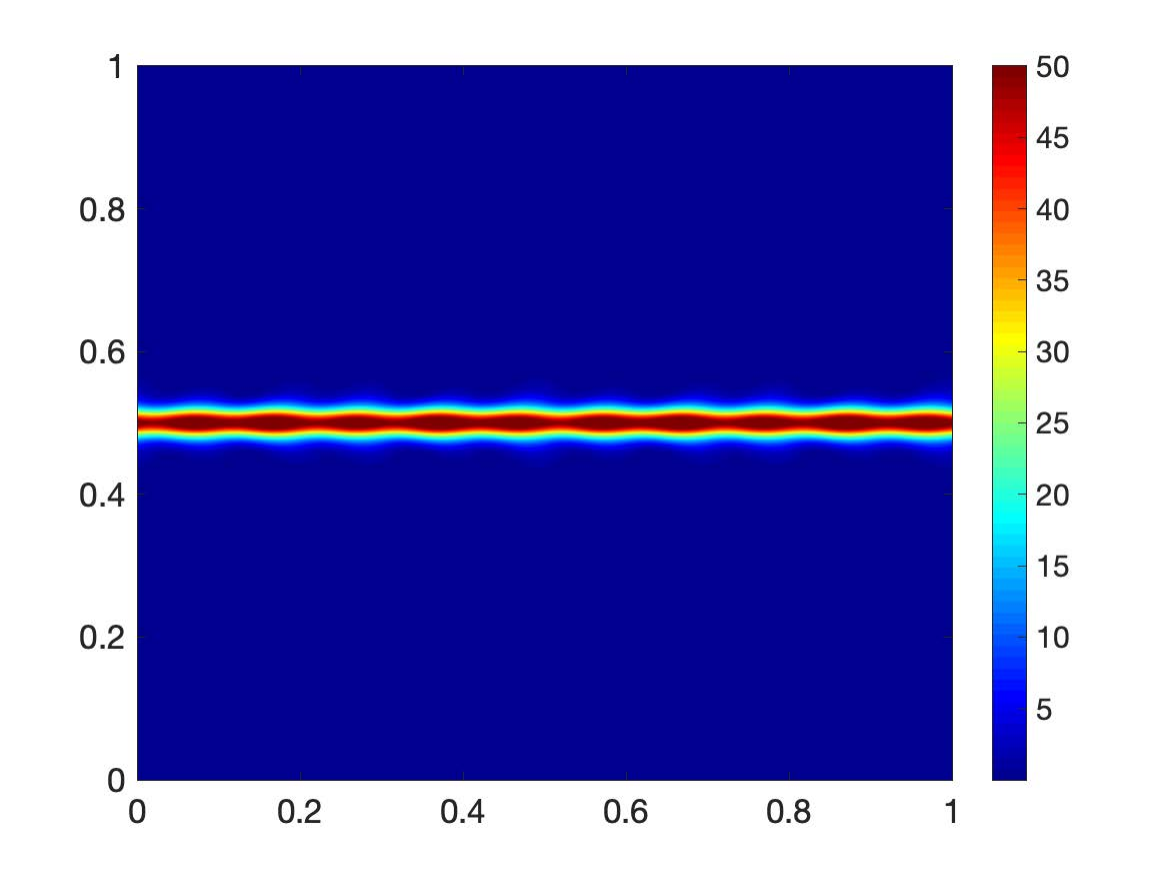}
		\includegraphics[width=0.3\textwidth,height=0.27\textwidth, trim=0 0 0 0, clip]{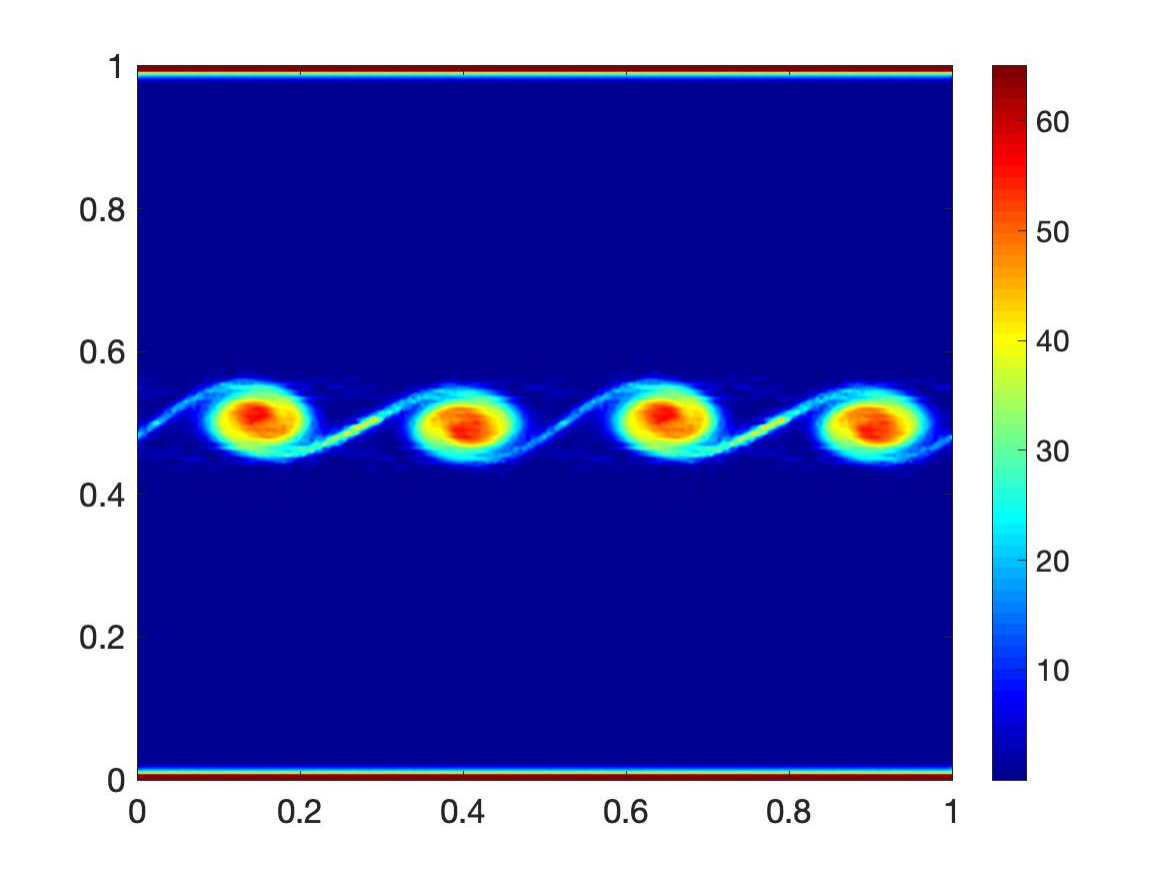}
		\includegraphics[width=0.3\textwidth,height=0.27\textwidth, trim=0 0 0 0, clip]{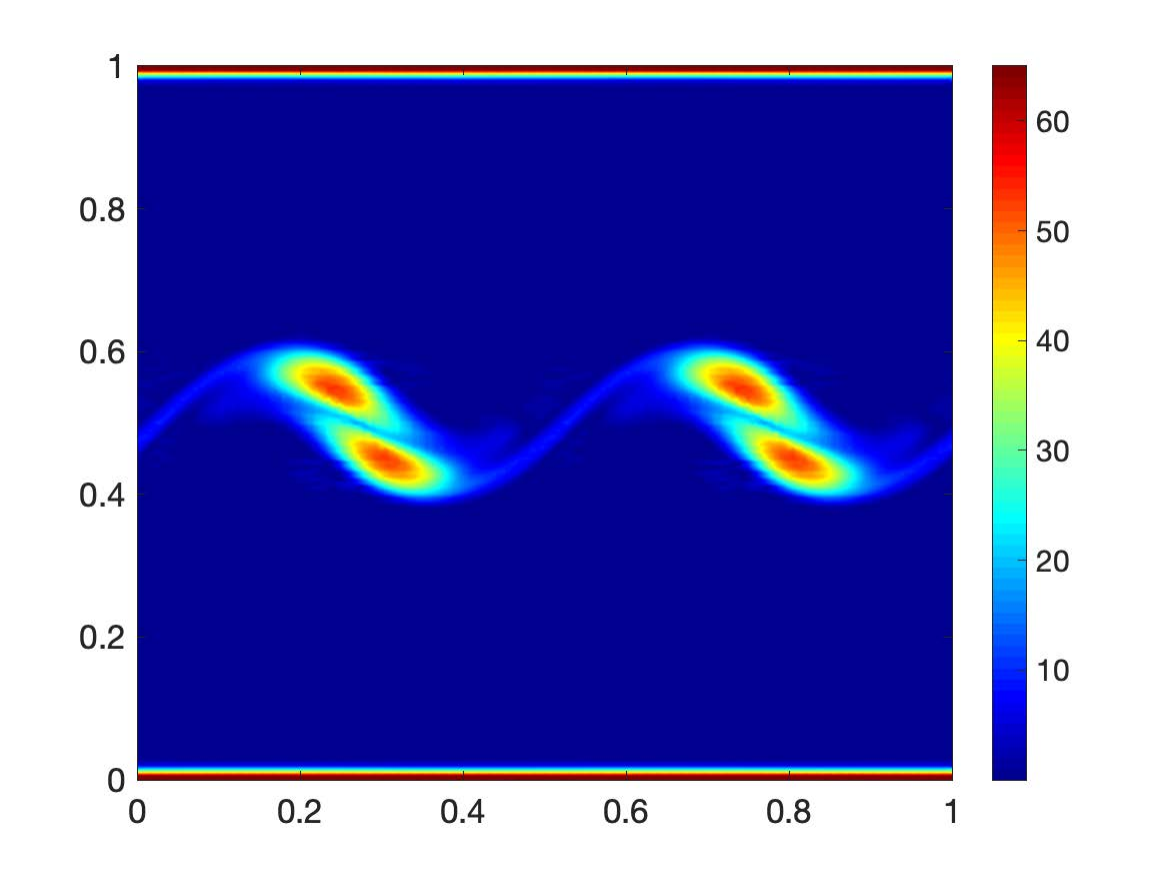}\\
				\includegraphics[width=0.3\textwidth,height=0.27\textwidth, trim=0 0 0 0, clip]{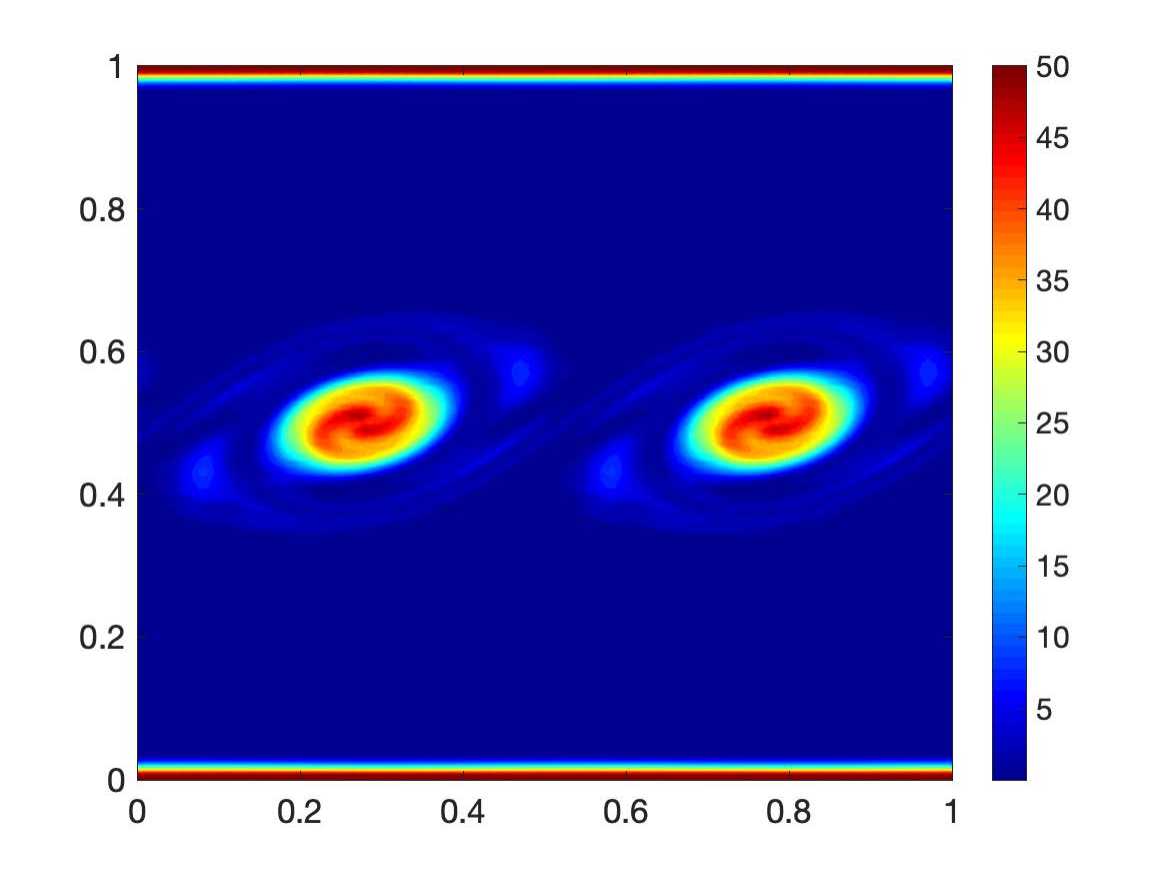}
		\includegraphics[width=0.3\textwidth,height=0.27\textwidth, trim=0 0 0 0, clip]{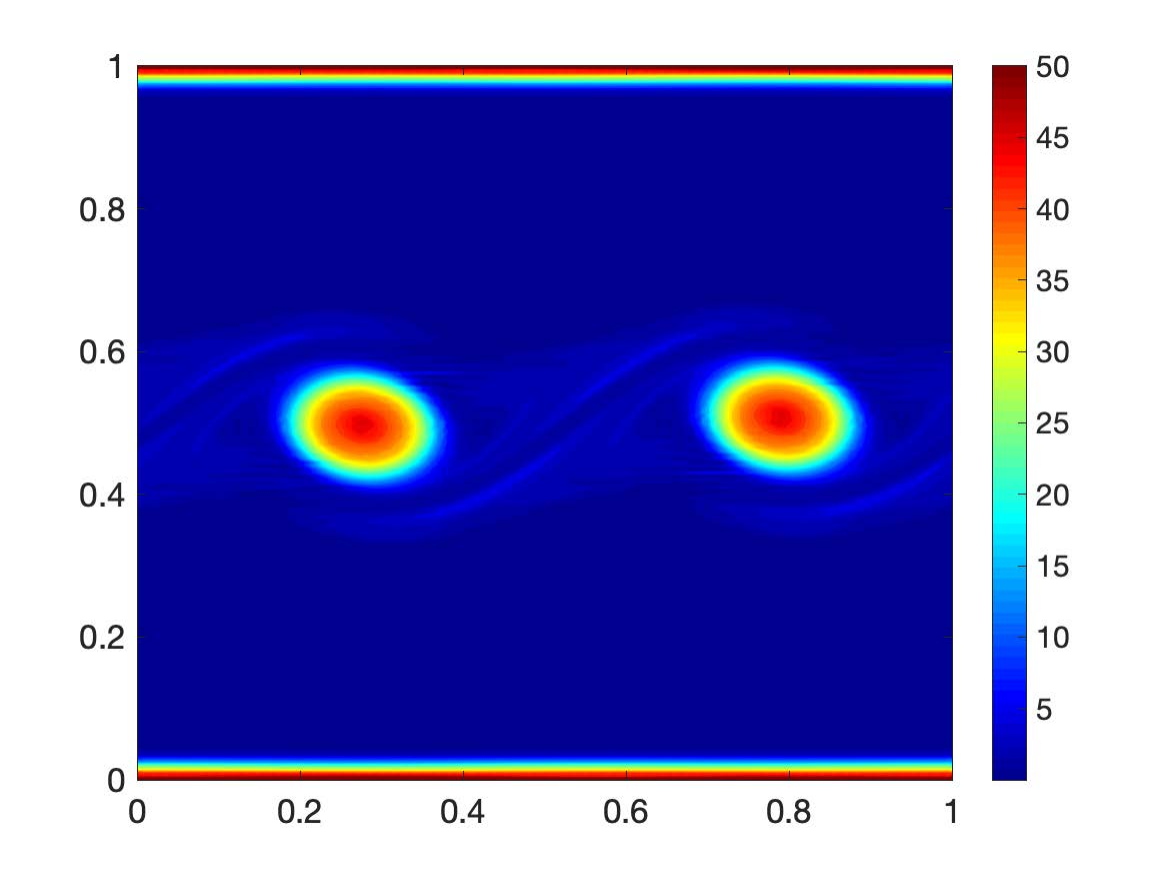}
		\includegraphics[width=0.3\textwidth,height=0.27\textwidth, trim=0 0 0 0, clip]{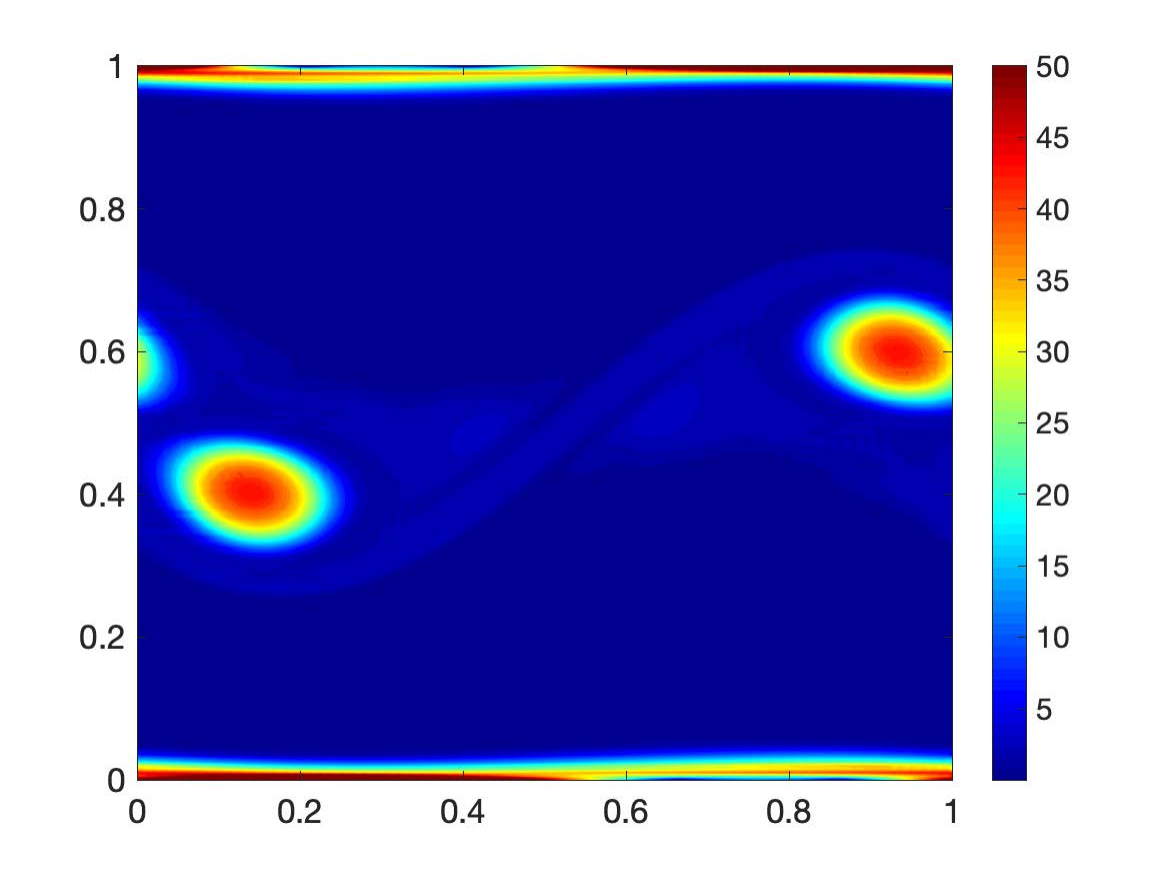}\\
				\includegraphics[width=0.3\textwidth,height=0.27\textwidth, trim=0 0 0 0, clip]{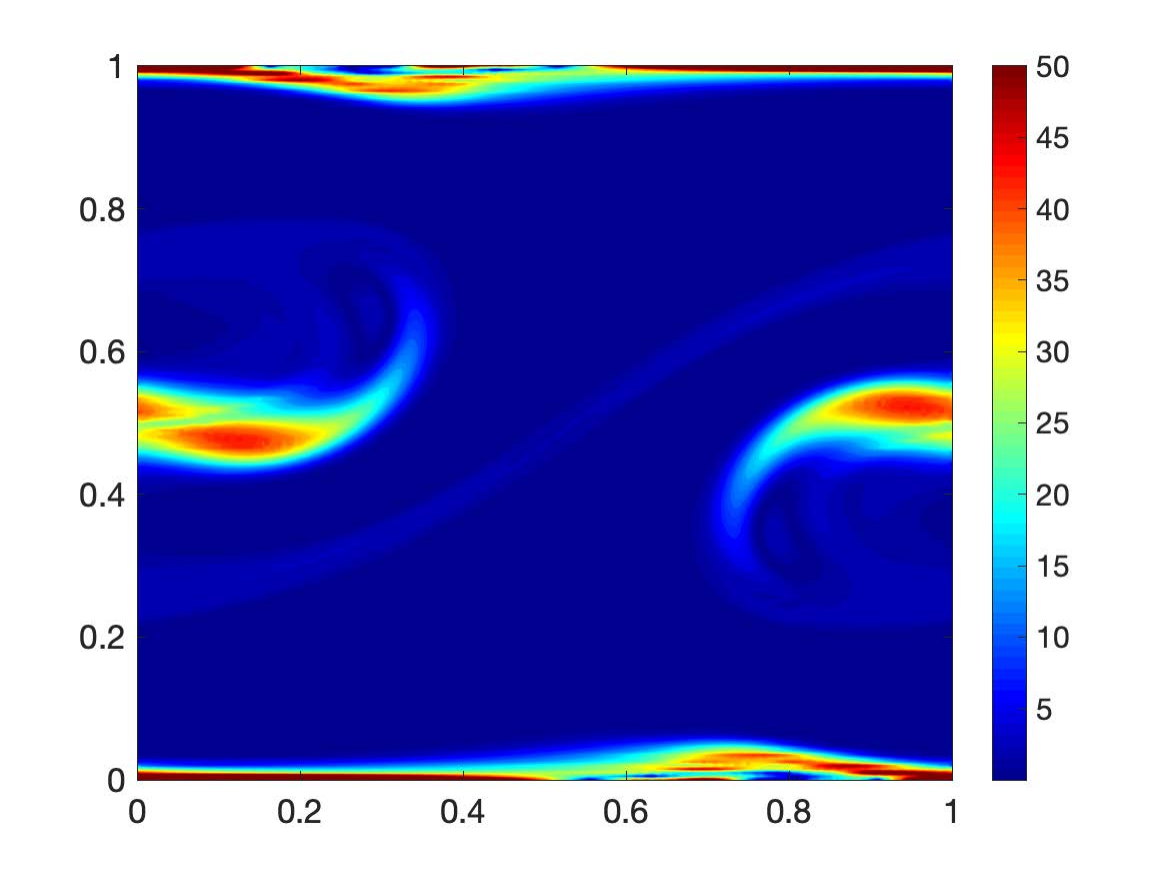}
		\includegraphics[width=0.3\textwidth,height=0.27\textwidth, trim=0 0 0 0, clip]{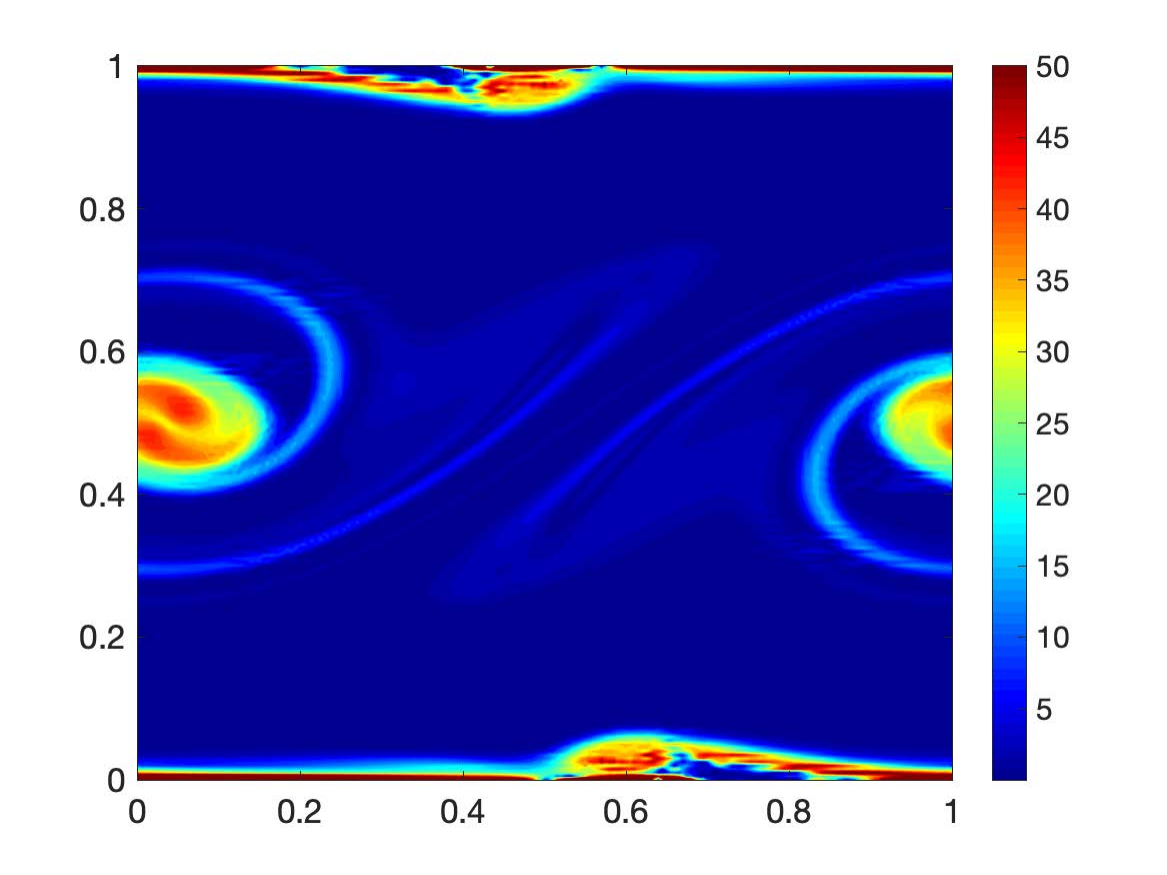}		
	\end{center}
	\caption{Shown above are the absolute vorticity contours of the solution velocity at $t=$0, 0.5, 1, 2, 3, 4, 4.5 and 5 (left to right, top to bottom). \label{KHsoln} }
\end{figure}

For our last test we consider a test problem from \cite{SJLLLS18} for simulating 2D Kelvin-Helmholtz instability.  The domain is the unit square, with periodic
boundary conditions at $x=0,1$, representing an infinite extension in the horizontal direction.  At $y=0,1$, we enforce for $t>0$ a no slip condition, which differs from  \cite{SJLLLS18} as they use a no penetration and free slip condition.  However, as these boundaries are far from the physical behavior of interest, there is little effect on the qualitative behavior of the solution.  The initial condition is set by
\[
\bu_0(x,y) = \left( \begin{array}{c} u_{\infty} \tanh\left( \frac{2y-1}{\delta_0} \right) \\ 0 \end{array} \right) + c_n \left( \begin{array}{c} \partial_y \psi(x,y) \\ -\partial_x \psi(x,y) \end{array} \right),
\]
where $\delta_0=\frac{1}{28}$ is the initial vorticity thickness, $u_{\infty}=1$ is a reference velocity, $c_n$ is a noise/scaling factor taken to be $10^{-3}$, and
\[
\psi(x,y) = u_{\infty} \exp \left( -\frac{(y-0.5)^2 }{\delta_0^2} \right) \left( \cos(8\pi x) + \cos(20\pi x) \right).
\]
The Reynolds number is defined by $Re=\frac{\delta_0 u_{\infty}}{\nu} = \frac{1}{28 \nu}$, and $\nu$ is defined by selecting $Re$.  We use $Re=$100 for our test.

We compute solutions for  EMAC discretized with $(P_2,P_1)$ Taylor-Hood elements on a $h=\frac{1}{128}$ uniform mesh, together with BDF2 time stepping and a time step size of $\Delta t=0.01$.  Solutions are computed up to $T=$5, with plots of vorticity contours shown in figure \ref{KHsoln} matching those in   \cite{SJLLLS18} qualitatively well.  

\begin{figure}[h]
	\begin{center}
		\includegraphics[width=0.49\textwidth,height=0.27\textwidth, trim=0 0 0 0, clip]{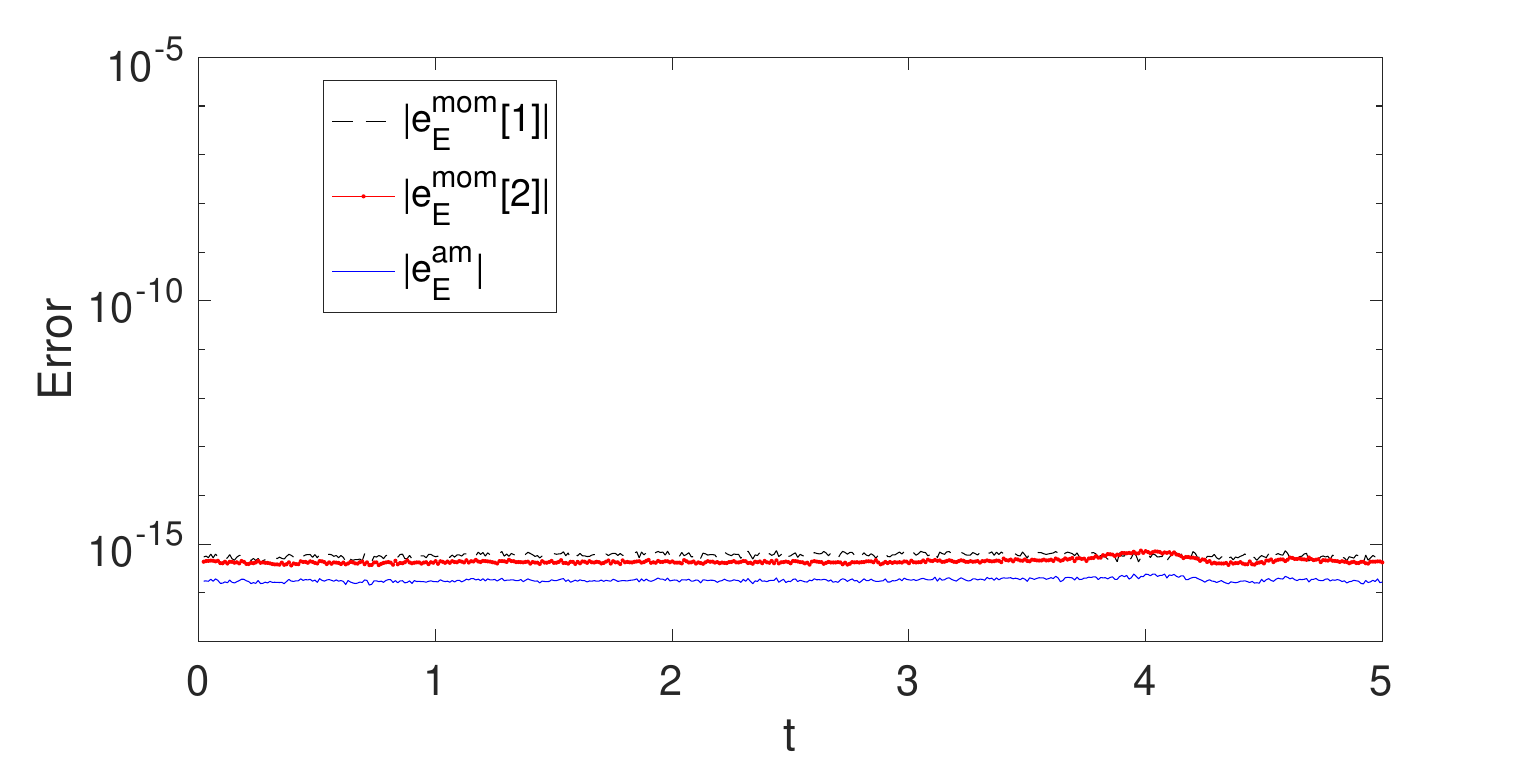}
		\includegraphics[width=0.49\textwidth,height=0.27\textwidth, trim=0 0 0 0, clip]{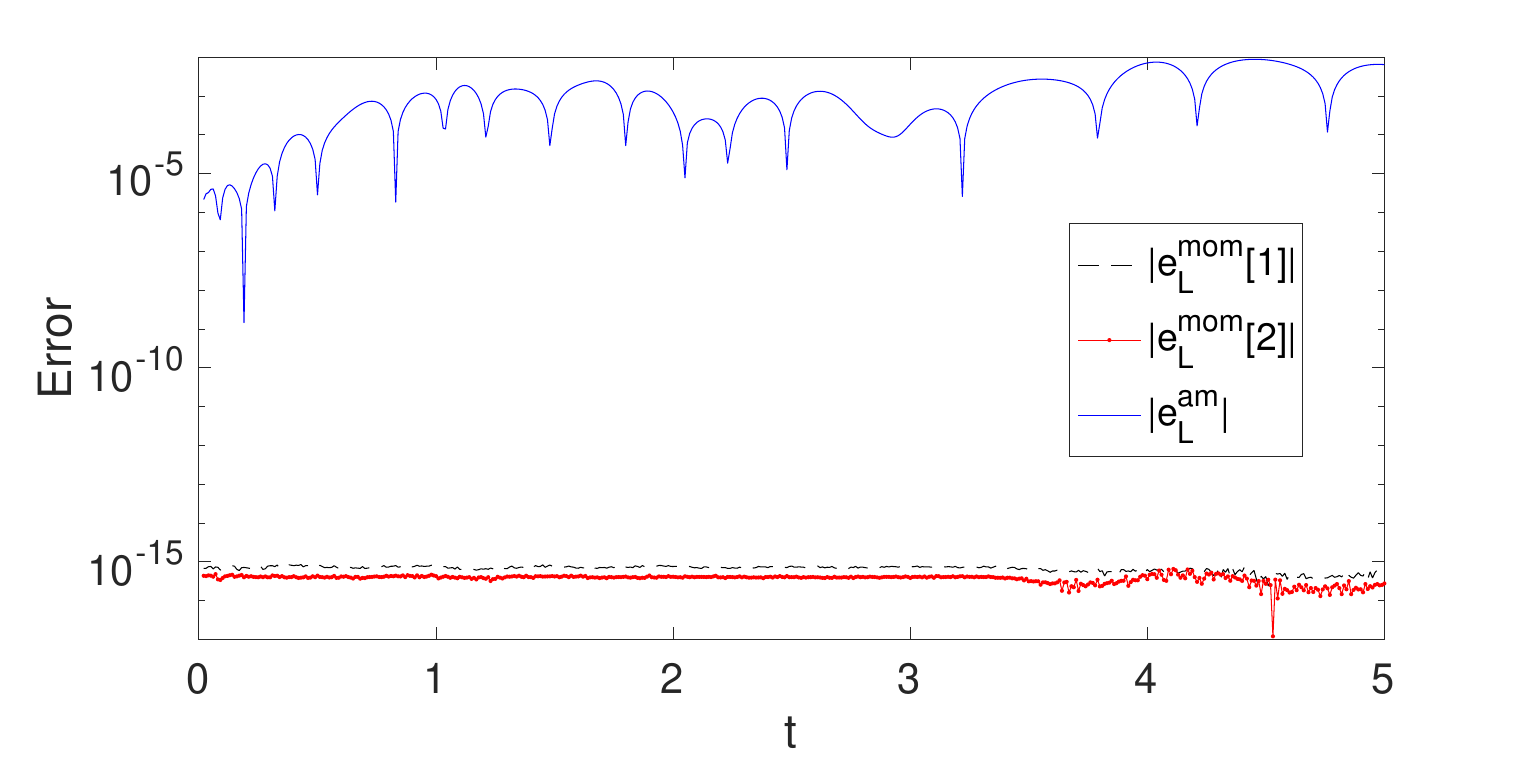}
	\end{center}
	\caption{Shown above is error in discrete local Eulerian (left) and Lagrangian (right) momentum and angular momentum conservation versus time in the Kelvin-Helmholtz problem.  \label{KHerr} }
\end{figure}

The subdomain $\omega$ is defined to be the square $[\frac18,\frac14] \times [\frac18,\frac14]$.  For this domain on this mesh, we have that $\omega_h=\omega$.  Plots of discrete Eulerian momentum and angular momentum are shown in figure \ref{KHerr} at the top, and we observe that these quantities are conserved exactly, just as the theory above predicts.  For discrete Lagrangian momentum and angular momentum, we solve the transport equation using backward Euler (BDF1), and plots of momentum and angular momentum are shown in figure \ref{KHerr} at bottom.  We observe exact local conservation of discrete Lagrangian momentum, and conservation of discrete Lagrangian angular momentum consistent with the discretization error, as predicted above.

%
%

\section{Future directions}

We have shown that continuous Galerkin discretizations of the Navier-Stokes equations using EMAC nonlinearity form admit (appropriately defined) exact local balances / conservation laws of momentum and angular momentum.  These discrete local balances are constructed as weak forms of the momentum and angular momentum conservation laws, and are equivalent to the usual conservation law definitions before discretization.  In the discrete case, however, these weak formulations are not equivalent to the usual conservation law definitions and even their Eulerian and Lagrangian constructions are not equivalent.  That the discrete schemes admit any exact local balances at all is quite rare, and we note that such an analysis is not possible for such common Navier-Stokes nonlinearity formulations such as convective, skew-symmetric or rotational.  We remark that the `conservative' formulation of nonlinear terms (referred to as CONS, utilizing $\Div(\bu_h \bu_h^T)$) also maintains the same conservation properties for momenta as EMAC. However, CONS fails to achieve a proper global energy balance, when $\Div\bu_h\neq0$, unlike EMAC. This deficiency leads to  unstable finite element schemes using CONS; see examples of CONS underperformance in \cite{CHOR17,CHOR19}.

Future directions for this work could include an extension of these ideas to other conservation laws of Navier-Stokes such as energy, helicity, enstrophy in 2D, vorticity, and others.  It is currently unclear to the authors how to construct appropriate local balances for these quantities for the continuous Galerkin method.

\section{Declaration of competing interest}

The authors declare the following financial interests/personal relationships which may be considered as potential competing interests: Both authors report financial support was provided by National Science Foundation.

\section{Acknowledgements}
M.O. was partially supported by the NSF grant DMS-2309197.
L.R. was partially supported by the NSF grant DMS-2152623.
We would like to express our gratitude to Thomas J.R. Hughes for his insightful discussions that inspired our research on the local conservation properties of EMAC.

\bibliographystyle{siam}
\bibliography{literature,bib_file}

\end{document}